\input amstex
\documentstyle{amsppt}
\magnification=\magstep1

\NoBlackBoxes
\TagsAsMath

\pagewidth{6.5truein}
\pageheight{9.0truein}

\long\def\ignore#1\endignore{#1}

\ignore
\input xy \xyoption{matrix} \xyoption{arrow}
          \xyoption{curve}  \xyoption{frame}
\def\edge{\ar@{-}}
\def\dttdar{\ar@{.>}}
\def\drbl{\save+<0ex,-2ex> \drop{\bullet} \restore}
\def\levelpool#1{\save [0,0]+(-3,3);[0,#1]+(3,-3)
                  **\frm<10pt>{.}\restore}
\def\dashedge{\ar@{--}}
\def\dropdown#1#2{\save+<0ex,-#1ex> \drop{#2} \restore}
\endignore

\def\seq{\mathrel{\widehat{=}}}
\def\la{{\Lambda}}
\def\lamod{\Lambda\text{-}\roman{mod}}
\def \len{\operatorname{length}} 
\def\AA{{\Bbb A}}
\def\PP{{\Bbb P}}
\def\CC{{\Bbb C}}
\def\SS{{\Bbb S}}
\def\ZZ{{\Bbb Z}}
\def\NN{{\Bbb N}}
\def\hom{\operatorname{Hom}}
\def\kleinhom{\operatorname{hom}}

\def\soc{\operatorname{soc}}

\def\stab{\operatorname{Stab}}

\def\maxtopdeg{\operatorname{max-topdeg}}

\def\GL{\operatorname{GL}}

\def\Ker{\operatorname{Ker}}
\def\Im{\operatorname{Im}}

\def\Schu{\operatorname{Schu}}

\def\Ext{\operatorname{Ext}}
\def\kleinext{\operatorname{ext}}

\def\A{{\Cal A}}
\def\frakA{{\frak A}}
\def\B{{\Cal B}} 
\def\C{{\Cal C}}
\def\frakC{{\frak C}}
\def\D{{\Cal D}}
\def\E{{\Cal E}}
\def\F{{\Cal F}}
\def\Gr{{\Cal G}r}

\def\I{{\Cal I}}

\def\T{{\Cal T}}
\def\L{{\Cal L}}

\def\R{{\Cal R}}
\def\S{{\sigma}}
\def\m{{\frak m}}
\def\s{{\frak s}}
\def\t{{\frak t}}
\def\T{{\Cal T}}

\def\aut{\operatorname{Aut}}

\def\autlap{\operatorname{Aut}_\la(P)}
\def\unirad{\bigl(\operatorname{Aut}_\la(P) \bigr){}_u}
\def\End{\operatorname{End}}

\def\U{{\Cal U}}

\def\dim{\operatorname{dim}}

\def\aux{{\operatorname{Schu}}}

\def\modlad{\operatorname{\bold{Mod}}^\Lambda_d}
\def\toptd{\operatorname{\bold{Mod}}^T_d}
\def\layS{\operatorname{\bold{Mod}}(\SS)}
\def\Hom{\operatorname{Hom}}
\def\grasstd{\operatorname{\frak{Grass}}^T_d}
\def\flagtd{\operatorname{\frak{Flag}}^T_d}
\def\autlap{\aut_\Lambda(P)}
\def\grassS{\operatorname{\frak{Grass}}(\S)}
\def\grassSS{\operatorname{\frak{Grass}}(\SS)}
\def\grass#1{\operatorname{\frak{Grass}}(#1)}

\def\degen{\le_{\deg}}

\def\id{\operatorname{id}}

\def\BongAdv{{\bf 1}}
\def\Bong{{\bf 2}}
\def\Bongtrond{{\bf 3}}
\def\GeomII{{\bf 4}}
\def\GeomIV{{\bf 5}}
\def\Bor{{\bf 6}}
\def\constructing{{\bf 7}}
\def\GabI{{\bf 8}}
\def\GabII{{\bf 9}}
\def\Har{{\bf 10}}
\def\dom{{\bf 11}}
\def\menace{{\bf 12}}
\def\class{{\bf 13}}
\def\Hum{{\bf 14}}
\def\Jou{{\bf 15}}
\def\KacI{{\bf 16}}
\def\KacII{{\bf 17}}
\def\Ki{{\bf 18}}
\def\Kol{{\bf 19}}
\def\Kraft{{\bf 20}}
\def\Kra{{\bf 21}}
\def\Rie{{\bf 22}}
\def\RieScho{{\bf 23}}
\def\Ros{{\bf 24}}
\def\SkoZwa{{\bf 25}}  
\def\Zwa{{\bf 26}}
\def\ZwaII{{\bf 27}}

\topmatter

\title Top-stable degenerations of finite dimensional
representations I
\endtitle

\rightheadtext{Top-stable degenerations}

\author  Birge Huisgen-Zimmermann \endauthor

\address Department of Mathematics, University of California, Santa
Barbara, CA 93106 \endaddress

\email birge\@math.ucsb.edu \endemail

\thanks The author was partly supported by a grant from the National
Science Foundation. 
\endthanks

\dedicatory  Dedicated to Claus Michael Ringel on the occasion of his
sixtieth birthday \enddedicatory

\abstract  Given a finite dimensional representation $M$ of a finite
dimensional algebra, two hierarchies of degenerations of
$M$ are analyzed in the context of their natural orders:  the poset of
those degenerations of
$M$ which share the top $M/JM$ with $M$  --  here $J$ denotes the
radical of the algebra  --  and the sub-poset of those 
which share the full radical layering $\bigl(J^lM/J^{l+1}M\bigr)_{l \ge 0}$
with $M$.  In particular, the article addresses existence of proper
top-stable or layer-stable degenerations  --  more generally, it addresses
the sizes of the corresponding posets including bounds on the lengths of
saturated chains  --   as well as  structure and classification. 
\endabstract

\endtopmatter

\document

\head 1. Introduction \endhead

Let $\la$ be a basic finite dimensional algebra over an algebraically
closed field $K$. Our objective is to
understand major portions of the degeneration theory of the finite
dimensional representations of $\la$.  This line of inquiry was triggered
by Gabriel's and Kac's initial exploitation of the affine variety
$\modlad$ that parametrizes the
$d$-dimensional representations of $\la$ (\cite{\GabI},
\cite{\GabII},
\cite{\KacI}, \cite{\KacII}) towards a better understanding of hereditary
algebras  --  equivalently, a better understanding of quiver
representations.  It was further pursued by Kraft, Riedtmann, Bongartz,
Schofield, Skowro\'   nski, Zwara, among many others (a selection of basic
references is
\cite{\Kraft}, \cite{\Kra}, \cite{\Rie}, \cite{\RieScho},
\cite{\BongAdv}, \cite{\Bongtrond}, \cite{\SkoZwa}, \cite{\Zwa}
\cite{\ZwaII}).   In intuitive terms, degenerations of a
$\la$-module $M$ result from a geometry-guided process of breaking its
`structural bonds', so as to incrementally unravel its
structure and relate it to less complex ones for the purpose of
analysis; for the definition, see Section 3.  Yet, in general, the problem
of classifying all degenerations of a given finite dimensional
representation $M$ is wild.  By contrast, there are strong indications that
the {\it top-stable\/} degenerations of $M$  (i\.e\., those degenerations
$M'$ which satisfy $M/JM \cong M'/JM'$, where $J$ denotes the Jacobson
radical of $\la$) should be classifiable, at least when
$T = M/JM$ is squarefree (that is, when $T$ has no simple summands of
multiplicity $\ge 2$).  In tandem, we address the smaller collection of all
{\it layer-stable\/} degenerations of $M$, namely those which share the
entire radical layering
$\SS(M) = (J^lM/ J^{l+1}M)_{l\ge 0}$ with $M$.  Our primary goal
in this article is to give a detailed general description of the
hierarchy of top-stable degenerations of $M$ in case
$T$ is squarefree (see Theorems 4.4, 4.5, and 5.1, 5.6) and to set up a
framework for classification (see Conjectures 4.6).  

Here is an excerpt of the basic theory:  Without loss of
generality, we may assume $\la$ to be a path algebra modulo relations,
$\la = KQ/I$, where $Q$  is a quiver and $I$ an admissible ideal in the
path algebra $KQ$. 
It turns out that, for any finite dimensional left $\la$-module $M$ with
squarefree top, the following two numerical invariants govern the size of
the poset of top-stable degenerations of $M$:  
\roster
\item"$\bullet$" The difference $\t - \s$, where $\t$ is the number of
simple summands in the top of $M$ and $\s$ the number of indecomposable
summands of $M$, and 
\item"$\bullet$" the difference 
$\m = \dim_K \Hom_\la(P,JM)- \dim_K \Hom_\la (M,JM)$, where $P$ is a
projective cover of $M$.
\endroster
We call a module {\it local\/} if it has a simple top.
The integer $\t-\s$, on one hand, measures the discrepancy from the case 
where $M$ is a direct sum
of local modules.  On the other hand, if $M \cong P/C$, the integer $\m$
measures the discrepancy from the situation where all homomorphisms $P
\rightarrow JP$ leave $C$ invariant.  The sum
$\m + \t - \s$ is a measure of how far $C$ is from being stable under
arbitrary endomorphisms of $P$, yet the two individual invariants, $\m$
and $\t -
\s$, play fundamentally different roles relative to the degeneration
theory of $M$.  We note that $\m$ is always bounded from
above by the sum of the multiplicities of the simple summands of $T$ as
composition factors of $JM/\soc(JM)$.     

\proclaim{Theorem A. Top-stable degenerations} Suppose $T =
M/JM$ is a direct sum of $\t$ pairwise non-isomorphic simple $\la$-modules
and $P$ a projective cover of $T$.  Write $M$ in the form $M \cong P/C$
with $C \subseteq JP$.

{\rm (1)} The lengths of chains of proper top-stable degenerations of
$M$ are bounded above by 

\qquad $\m + \t - \s$.
\smallskip

{\rm (2)} Existence:  The following conditions are equivalent.
\roster

\item"(i)"  $M$ has a proper top-stable degeneration.

\item"(ii)"  $\m + \t - \s > 0$.

\item"(iii)"  Either $M$ is not a direct sum of local modules, or else
$C$ fails to be invariant under all homomorphisms $P \rightarrow JP$. 

\item"(iv)" $C$ fails to be invariant under all endomorphisms of $P$.
\endroster

{\rm (3)} Unique existence:  $M$ has a unique proper top-stable
degeneration if and only if $M$ is a direct sum of local modules and $\m =
1$.  

If $\m = 0$ and $\t - \s = 1$, $M$ has precisely two distinct proper
top-stable degenerations. {\rm (}For all values 
$\m + \t - \s \ge 2$, infinitely many top-stable degenerations can be
realized.{\rm )}
\smallskip

{\rm (4)} Bases:  If $M'$ is a top-stable degeneration of $M$, then $M$
and $M'$ share a basis consisting of paths in the quiver
$Q$.

To make this statement more precise, we assume {\rm (}without loss of
generality{\rm )} that $P$ is a left ideal of $\la$.  Then there
exist submodules $D$ and $D'$ of $P$ with $M \cong P/D$ and $M' \cong
P/D'$, together with a set
$\B$ of paths in $KQ$ {\rm{(}}which is closed under right subpaths{\rm{)}},
such that
$\{ p + D \mid p \in \B \}$ is a basis for $P/D$ and $\{ p + D' \mid p \in
\B \}$ is   a basis for $P/D'$.
\smallskip

{\rm (5)} The maximal top-stable degenerations of $M$ always possess a
fine moduli space, $\maxtopdeg(M)$, classifying them up to isomorphism. 
This moduli space is a projective variety of dimension at most
$\max\{0,\, \m + (\t -\s) -1\}$.
\smallskip

{\rm (6)} The case $\m = 0$:
Let $M = \bigoplus_{1 \le k \le \s} M_k$ be a
decomposition of $M$ into indecomposable summands.  Then $M$ has only
finitely many top-stable degenerations, and each of them is a direct sum of
top-stable degenerations of the $M_k$.  Moreover, $M$
degenerates top-stably to $M'$ precisely when
$M \le M'$ in the $\Ext$-order  {\rm{(}}see beginning of
Section {\rm{3)}}. 

In particular: If $\hom_\la(P, JP) = \hom_\la(P,\soc(JP))$, the
degeneration order on the class of left $\la$-modules with top $T$
coincides with the $\Ext$-order.
\endproclaim  

The upper bounds in Theorem A(1),(5) are attained for arbitrary values
of
$\m$ and $\t - \s$. 

By way of preview of work in progress, we mention that,  in the situation
of repeated simple summands in the top $T$, the pressure placed on the
structure of $M$ by non-existence of proper top-stable
degenerations exceeds expectations based on part A(2) of the
theorem:  In fact, suppose that
$T = M/JM \cong  S_1^{t_1} \oplus \cdots \oplus S_n^{t_n}$, where the $S_i$
are pairwise non-isomorphic simples $\la e_i/Je_i$ for suitable primitive
idempotents $e_i$, and let
$P$ again be a projective cover of $T$.  Then
$M$ is devoid of proper top-stable degenerations if and only if the
following three conditions are satisfied:

$\bullet$  $M$ is a direct sum of local modules, say $M = \bigoplus_{1 \le
i \le n} \bigoplus_{1 \le j \le t_i} M_{ij}$, where $M_{ij} = \la
e_i/C_{ij}$, i.e., $\bigl( \sum_{1 \le i \le n} t_i \bigr) - \s = 0$,

$\bullet$ $\dim \hom_\la(P, JM) = \dim \hom_\la(M,JM)$, i\.e\., $\m = 0$,
and

$\bullet$  For each $i \le n$, the $C_{ij}$ are linearly ordered with
respect to inclusion.
\medskip 

\proclaim{Theorem B. Layer-stable degenerations} Let $M$, $\s$, and $\t$ be
as in Theorem A.

{\rm (1)} If $M$ is a direct sum of local modules, i\.e., if
$\t - \s = 0$, then $M$ has no proper layer-stable degenerations.
\smallskip

{\rm (2)}  If $\m = 0$, then every minimal layer-stable
degeneration of $M$ {\rm (}see {\rm 3.7)} is of the form $M = U \oplus
M/U$, where $U \subseteq M$ is a layer-stably embedded submodule
{\rm (}meaning $J^l U = J^l M \cap U$ for all $l${\rm )}.  
\smallskip

{\rm (3)}  On the other hand, for any positive integer $r$, there exists
an indecomposable finite dimensional module $M$ with squarefree top, over a
suitable finite dimensional algebra $\la$, such that $M$ has a
$\PP^r$-family of pairwise non-isomorphic indecomposable layer-stable
degenerations.
\endproclaim

Theorems A and B are excerpted from Theorems 5.1 and 5.6. For the case of
local modules (i.e., $\t=1$), see Theorem 4.4. The second major result of
Section 4, Theorem 4.5, addresses the vertical
and horizontal size of the poset of top-stable degenerations of a local
module.  In rough terms, it says that ``all" phenomena not excluded by
either Theorem A or the conjectures outlined below actually occur. 
In particular, this means that every finite direct product of projective
spaces can be realized as the moduli space, $\maxtopdeg(M)$, classifying
the maximal top-stable degenerations of a local module $M$.  Moreover,
severe deviation from catenarity is a common phenomenon.  The conjectures
are backed by strong experimental evidence,  collected with the
computational methods which we briefly sketch prior to the examples in
Sections 4.B and 5.B.  These methods will be further elaborated in joint
work with R\. Donagi
\cite{\constructing}.
\medskip

{\bf Conjectures.}  (Excerpt from Conjectures 4.6) We keep the notation of
Theorem A.  If
$M = P/C$ has squarefree top of length $\t$, the following hold:

(1)  There exist finitely many non-negative intergers $n_1 \le 
\cdots \le n_s$ such that conditions (a) and (b) below are satisfied. 

(a)  The irreducible components of the moduli space $\maxtopdeg(M)$
coincide with its connected components, and each component is
isomorphic to a direct product of projective spaces, $\prod_{i \in I}
\PP^{n_i}$, for a suitable subset $I$ of $\{1, \dots, s\}$.

(b)  For each $h \le m$, the top-stable degenerations of height $h$ above
$M$ (in the usual degeneration order) can be partitioned into finitely
many disjoint families, each of which is bijectively parametrized by an
orbit of the standard torus action on a product $\prod_{i \in I}
\PP^{n_i}$, for some $I \subseteq \{1, \dots, s\}$.  (The torus we
refer to is $\prod_{i \in I} (K^*)^{n_i + 1}$.)

If $s$ and $\sum_{i \le s} n_i$ are chosen minimal with respect to (a) and
(b), then the numbers $s$ and $n_1, \dots, n_s$ are isomorphism invariants of
$M$. 
\smallskip

(2)  Suppose, more restrictively, that $M$ is local (in particular, this
implies $\t -\s = 0$).  Then
$\maxtopdeg(M)$ is irreducible.  If, in addition, we require the poset of
top-stable degenerations of $M$ to be finite, its cardinality (counting
$M$) is bounded from below by $\m +1$ and from above by $2^{\m}$. 
Moreover, $M$ has a unique maximal top-stable degeneration in that case,
and all saturated chains linking $M$ to the latter have length $\m$. 
\smallskip 

Our examples will show that the catenarity condition under (2) fails
without the finiteness hypothesis  --  with arbitrarily large gaps between
the lengths of maximal and minimal saturated chains in
general.  Moreover, all of the conjectural bounds are attained. 
\medskip   

An ulterior motivation for the study of
degenerations stems from their impact on possibilities of more
comprehensive classification, in particular, on the problems of when the
full classes of modules with fixed top $T$ or fixed radical layering $\SS$
admit fine or coarse moduli spaces (see
\cite{\class}).  In fact, as is well-known, degenerations of modules
constitute the foremost obstruction in the way of helpful quotients of
stable subvarieties of
$\modlad$ by their $\GL_d$-actions, for the following reason:  As the
$\GL_d$-orbits are in one-to-one correspondence with the isomorphism types
of $d$-dimensional $\la$-modules, one looks for quotients that separate
orbits; such separation is obviously ruled out if the orbits fail to be
closed, that is, if the modules represented by these orbits have proper
degenerations.  Unfortunately  --  in contrast to the
$G$-spaces typically arising in classical geometric problems  --  there are
hardly any closed orbits in the module-theoretic setting, and
hence the standard methods of invariant theory are not effective in this
context.  This was already addressed by King (see \cite{\Ki}). 

Concerning methodology, the class of top-stable degenerations of a given
module admits a twofold approach:  one via the quasi-affine subvariety
$\toptd$ of the classical variety $\modlad$ that para\-metriz\-es the
$d$-dimensional $\la$-modules with fixed top $T$, the other -- 
our main resource here  --  via a projective variety $\grasstd$ which
parametrizes this same class of representations.  The variety $\grasstd$
is in turn partitioned by a group action whose orbits correspond
bijectively to the isomorphism classes of the $d$-dimensional
representations with top $T$.   But, in contrast to the
reductive group acting on
$\modlad$, the group acting on $\grasstd$ has a large unipotent
radical in all interesting cases.  This fact, in combination with
completeness of $\grasstd$, brings a different cache of geometric methods
to bear, as witnessed for instance by Proposition 4.2 and Corollary 5.4,
next to the proofs of the main results.  The two modes of approach combine
smoothly, due to a result of Bongartz and the author (see Section 2).  

In Section 2 we summarize prerequisites; the more technical ones will
only be required for proofs and examples.  In Section 3 we assemble a
number of general facts about degenerations, not restricted to modules with
squarefree tops, of which 3.4--3.7 appear to hold the highest independent
interest.  Crucial definitions concerning degenerations, in particular the
concept of {\it height\/}, can be found at the end of Section 3.  Section
4 is focused on the local case, i\.e., the case of a simple top
$T$, a situation which stands out through numerous special features; since
all saturated chains of top-stable degenerations pass through direct sums
of local modules, the ``final degeneration stages" of nonlocal
modules lead back to the local scenario.  Section 5 addresses arbitrary
squarefree tops.  In Sections 4.B and 5.B we supplement the theory with
examples. 
\medskip

\noindent{\it Acknowledgements:\/} This research was partially supported by
a grant from the National Science Foundation.  The author wishes to thank Ron
Donagi for contributing Lemma 4.3 and J\'anos Koll\'ar for
Proposition 3.6.  Moreover, she is indebted to
Frauke Bleher and Sverre Smal\o\ for their many helpful comments on a
preliminary draft of this paper.  Finally, her thanks go to the referee for
the most thorough report she has ever received, containing several helpful
suggestions and alerting her to some gaps in the (original) arguments.

\head 2. Prerequisites \endhead 

\noindent{\it General conventions:\/} 

Throughout, $\la$ will be a basic finite dimensional algebra over an
algebraically closed field $K$.  Hence, we may assume without loss of
generality that
$\la = KQ/I$, where $Q$ is a quiver and
$I$ an admissible ideal in the path algebra $KQ$.  Moreover, $L$ will 
denote the largest integer for which the power
$J^L$ of the Jacobson radical $J$ of $\la$ does not vanish; in other words,
$L+1$ is the Loewy length of
$\la$.  The quiver $Q$ provides us with a distinguished set of primitive
idempotents
$e_1, \dots, e_n$ of
$\la$, which are in bijective correspondence with the vertices of
$Q$; notationally, we will often not distiguish between the vertices and
the $e_i$.  As is well-known, the factors $S_i =
\la e_i /J e_i$ then form a set of representatives for the simple left
$\la$-modules.   An element
$x$ of a (left) $\la$-module $M$ will be called a {\it top element\/} of
$M$ if $x \notin JM$ and $x$ is normed by some $e_i$, that is, $x = e_i
x$. 

Moreover, we will observe the following conventions:  The product $pq$ of
two paths $p$ and $q$ in $KQ$  stands for ``first $q$, then $p$"; in
particular, $pq$ is zero unless the end point of $q$ coincides with the
starting point of $p$.  In this spirit, we call a path
$p_1$ a {\it right subpath\/} of $p$ if $p = p_2 p_1$ for some path
$p_2$.  A {\it path in
$\la$\/} will be an element of the form $p +I$, where
$p$ is a path in $K Q$; we will suppress the residue notation, provided
there is no risk of ambiguity.  Further, we will gloss over the
distinction between the left
$\la$-structure of a module $M \in \lamod$ and its induced
$KQ$-structure.  (On the other hand, when we refer to paths in
$KQ$, we will always point this out explicitly.) 

Let $a_1, \dots, a_r$ be a set of algebra generators for
$\la$ over $K$; for $\la = KQ/I$ as above, a typical set
of generators consists of the primitive idempotents $e_1, \dots, e_n$
together with the (residue classes in $\la$ of the) arrows in
$Q$.  Recall that, for $d \in \NN$, the classical affine variety of
$d$-dimensional representations of
$\la$ can be described in the form
$$\modlad = \{(x_i)_{i \le r} \in \prod_{1 \le i \le r} \End_K(K^d)
\mid \text{\ the\ } x_i \text{\ satisfy all relations satisfied by the\ }
a_i\}.$$  
Given a semisimple left $\la$-module $T$, we denote by $\toptd$ 
(or $\operatorname{\bold{Mod}}^{\la,T}_d$) the locally closed subvariety of
$\modlad$ which contains precisely the points representing modules $M$
with top
$T$ (that is, with $M/JM \cong T$).  Thus the orbits of the
$\GL_d$-conjugation action on
$\toptd$ are in one-to-one correspondence with the isomorphism types of
$d$-dimensional left $\la$-modules that have top $T$. 
\bigskip

\noindent{\it Background on Grassmannians of $d$-dimensional
modules with fixed top $T$:\/} 

We give a brief review of those constructions and results from 
\cite{\class} which will be pivotal in the sequel (see also
\cite{\GeomII} and \cite{\GeomIV} for background).  For that purpose, we
fix a finite dimensional semisimple module $T$, together with a positive
integer $d$, and let $P$ be a projective cover of $T$.  

Our principal tool will be the
following closed subvariety of the classical Grassmannian $\Gr(d',JP)$ of 
$d'$-dimensional subspaces of the $K$-space $JP$, where $d' = \dim_K P
- d$:
$$\grasstd = \{C \in \Gr(d',JP) \mid C \text{\ is a\ }
\la\text{-submodule of\ } JP \}.$$ 
This variety comes with an obvious surjection
$$\phi: \grasstd \longrightarrow \{\text{isomorphism types of\ } d\text
{-dimensional modules with top\ } T\},$$  
sending $C$ to the class of $P/C$.  Clearly, the fibres of $\phi$ coincide
with the orbits of the natural $\autlap$-action on
$\grasstd$.   While the global geometry of the projective variety
$\grasstd$ cannot be reasonably compared with that of the quasi-affine
variety $\toptd$, the ``relative geometry'' of the
$\autlap$-stable subsets of $\grasstd$ is tightly linked to that of the
$\GL_d$-stable subsets of $\toptd$ (see
\cite{\GeomIV}).  Namely, if the orbits $\autlap.C$ and $\GL_d.x$
in $\grasstd$ and $\toptd$, respectively, represent the same
$\la$-module up to isomorphism, the assignment  $\autlap.C$
$\mapsto$ $\GL_d.x$ induces a bijection
$$\{ \autlap\text{-stable subsets of\ } \grasstd \} \rightarrow
\{\GL_d\text{-stable subsets of\ } \toptd \}$$  which preserves openness,
closures, connectedness, irreducibility, and types of singularities.  
This correspondence permits transfer of information between the two
settings. 
     
The first of the two isomorphism invariants of modules we use to organize
$\grasstd$ from a represen\-ta\-tion-theoretic viewpoint is the {\it
radical layering\/}.  It leads to a
partition of
$\grasstd$ into locally closed subvarieties (in general, it is not a
stratification, however).  

Recall that the {\it radical layering\/} of a module $M$ is the sequence
$$\SS(M) = (M/JM, JM/J^2M, \dots, J^{L-1}M/J^L M, J^L M),$$   
where we
identify semisimple modules with their isomorphism classes; in effect, we
are dealing with a matrix of discrete invariants keeping count of
the multiplicities of the simple modules in the individual radical layers
of $M$.  Correspondingly, we consider the following locally closed
subvarieties of
$\grasstd$:  
$$\grassSS  = \{ C \in \grasstd \mid \SS(P/C) = \SS\}$$ 
for any
$d$-{\it dimensional semisimple sequence $\SS$ with top\/}
$T$; by this we mean any sequence  
$(\SS_0, \dots, \SS_L)$ of semisimple modules such that
$\SS_0 = T$, each $\SS_l$ embeds into $J^lP/J^{l+1}P$, and
$\sum_{0 \le l \le L} \dim \SS_l = d$.  Clearly, the above one-to-one
correspondence between the $\autlap$-stable subsets of
$\grasstd$ and the $\GL_d$-stable subsets of $\modlad$ restricts to a
correspondence between the $\autlap$-stable subsets of $\grassSS$
and the
$\GL_d$-stable subsets of a locally closed subvariety
of $\toptd \subseteq \modlad$, which we denote by $\layS$:  It consists
of those points in the classical variety $\toptd$ which correspond to the
modules with radical layering $\SS$.  
\bigskip

\noindent{\it Further prerequisites, required for proofs
and examples only\/}: 

For the remainder of this section, we will assume $T$ to be
{\it squarefree\/}.  Thus
$$T = \bigoplus_{1 \le i \le \t} S_i,$$
after potential reordering of the primitive idempotents.  Correspondingly,
we may take
$P$ to be $\bigoplus_{1 \le i \le \t} \la e_i$.  In other words, $P$ is the
left ideal of $\la$ generated by $e_1 + \dots + e_\t$.

The restriction to squarefree $T$ allows for a simplified description of
the second invariant of a module, which we use to organize $\grasstd$ for
our purposes, namely its set of {\it skeletons\/}.  Since each
such skeleton consists of $d$ paths of length at most $L$ in the quiver
$Q$, the set of all skeletons of $d$-dimensional modules is finite.    

\definition{Definitions} For $C \in \grasstd$, we call a set
$\sigma$ of paths of lengths $\le L$ in $KQ$  a {\it skeleton\/} of $M =
P/C$  if the following two conditions are satisfied:
\roster

\item"(a)" $\sigma$ is closed under right subpaths, that is `$p_2 p_1 \in
\S \implies p_1 \in \S$';

\item"(b)" for each $l \in \{0, \dots, L\}$, the set 
$$\{ p + C \mid p \in \sigma,\ \len(p) = l \}$$ induces a $K$-basis for
$J^l M/ J^{l+1}M$.
\endroster

 Existence of at least one skeleton of $M$ is
obvious.  Note, moreover, that
condition (b) is equivalent to postulating that, for {\it every\/}
sequence $x_1, \dots, x_\t$ of top elements of
$M$ with $e_i x_i = x_i$ and any $\l \in \{0, \dots, L\}$, the set 
$$\bigcup_{1 \le i \le \t} \{px_i \mid p = p e_i \in \sigma,\ \len(p) = l
\}$$  induces a basis for $J^l M/ J^{l+1}M$.  This guarantees that
isomorphic modules have the same skeletons.

We further observe that (a) and (b) entail the following condition (c):
$\sigma \cap \{e_1, \dots, e_n\} = \{e_1, \dots, e_\t\}$.  To obtain a
suitable partition of $\grasstd$, we refer to a set
$\sigma \subseteq KQ$  as an {\it abstract
$d$-dimensional skeleton\/} with top $T$ if
$\sigma$ consists of $d$ distinct paths in $KQ$ satisfying conditions
(a) and (c).  An abstract skeleton need not be
realizable as a skeleton of a module.    

Every abstract
$d$-dimensional skeleton $\sigma$ with top $T$ gives rise to an
$\autlap$-stable locally closed subset of $\grasstd$, namely 
$$\grassS : = \{ C \in \grasstd \mid \sigma \text{\ is a skeleton of\ }
P/C\}.$$  
A priori, we do not exclude those skeletons $\sigma$ for which $\grassS$
is empty.

Each $\grassSS$ is the union of certain $\grassS$. Given $\SS$, the
appropriate choice of skeletons $\sigma$ is as follows:  An abstract
skeleton
$\sigma$ is said to be {\it compatible with\/} $\SS$ in case, for each
$l \in \{0, \dots, L\}$ and $i \in \{1, \dots, n\}$, the cardinality of
the set of paths in $\sigma$ that have length $l$ and end in
$e_i$ equals the multiplicity of the simple module $S_i$ in $\SS_l$.  If
$\sigma$ is a skeleton of $M$, then clearly $\sigma$ is compatible with
$\SS(M)$, and hence $\grassSS$ is the union of those $\grassS$ which
correspond to skeletons compatible with $\SS$. 
\enddefinition

\proclaim{Theorem 2.1} \cite{\class, Theorem 3.5}  Let $\SS$ be a
$d$-dimensional semisimple sequence with squarefree top $T$. Then the
sets $\grassS$, where $\S$ runs through the skeletons compatible with
$\SS$, are $\autlap$-stable and form an affine open cover of the variety
$\grassSS$ representing the modules with radical layering $\SS$.
\qed
\endproclaim

The $\autlap$-stability of the resulting cover $\bigl( \grassS
\bigr)_{\sigma}$ of $\grasstd$ will be pivotal in the sequel.  Caveat: 
While they are open in the pertinent subvarieties $\grassSS$, the affine
charts $\grassS$ fail to be open in $\grasstd$ in general.  

\definition{Connection with Schubert cells} We briefly comment on how this
affine cover of $\grasstd$ relates to the Schubert cells of the classical
Grassmannian $\Gr(d', JP)$ in which $\grasstd$ is embedded as a closed
subvariety (recall that $d' = \dim_KP - d$).  Each
$\grassS$ embeds
into the big Schubert cell, $\Schu(\sigma)$, which consists of
the points in $\Gr(d',JP)$ that complement
$\bigoplus_{p \in \S} Kp$ in $P$.  In fact, $\grassS$ is the
intersection of $\Schu(\sigma)$ with $\grassSS$, where $\SS$ is the
unique semisimple sequence with which $\sigma$ is compatible.  (We point
out that, by contrast, the intersection of $\Schu(\sigma)$ with $\grasstd$
is usually not closed under the
$\autlap$-action.)  In essence, the affine coordinates for
$\grassS$ presented below are the Pl\"ucker coordinates relative to a
conveniently chosen basis for $JP$ that contains $\sigma$.  Yet, we only
need to retain a small portion of the numerical data recorded by these
latter coordinates, owing to the fact that we are exclusively interested in
sub{\it modules\/} of
$JP$.  This restriction to submodules is also responsible for the fact
that our coordinates are determined by polynomials that can readily be
derived from the relations which we factor out of $KQ$ to obtain $\la$.
\enddefinition    

We conclude the section by describing the mentioned coordinate system for
$\grassS$, relevant only for our analysis of examples.  Starting with the
quiver
$Q$ and a finite left ideal generating set $\R$ of $I$, we can obtain
polynomials describing
$\grassS$ as outlined below;  proofs can
be found in \cite{\class, 3.10--3.13}.  A
$\sigma$-{\it critical pair\/} is a pair
$(\alpha, p)$, where $p$ is a path in $\sigma$ and $\alpha$ is an arrow
such that $\alpha p$ is a path which fails to belong to $\S$.   Given $C
\in \grassS$, there are unique scalars $c_{\alpha p, q}$ with 
$$\alpha p + C = \sum_{q \in \sigma(\alpha,p)}  c_{\alpha p,q} (q + C),$$
where $\sigma(\alpha, p)$ is the set of all those paths $q \in \sigma$
which are at least as long as $\alpha p$ and have the same
end point as $\alpha p$.  Clearly, the sets $\S(\alpha,p)$ for distinct
critical pairs $(\alpha,p)$ will intersect in general.  

We let
$N$ be the {\it disjoint\/} union of the $\sigma(\alpha,p)$, where
$(\alpha,p)$ traces the $\sigma$-critical pairs.  Whenever this
promotes clarity, we will use the indexing $q = q(\alpha,p)$  to
achieve formal disjointness. 

Given $C \in \grassS$, the point
$c = (c_{\alpha p, q}) \in \AA^N$ clearly determines $C \subseteq P =
\bigoplus_{1
\le i \le \t} \la e_i$.  Indeed, 
$$C = \sum_{(\alpha, p)\ \sigma{\text{-critical}}} \la
\bigl(\alpha p\ -
\sum_{q \in \sigma(\alpha, p)} c_{\alpha p,q} q \bigr).$$ 
Due to
\cite{\class, Theorem 3.14}, the assignment $\grassS \rightarrow
\AA^N$ that sends any point $C \in \grassS$ to the corresponding point $c
\in \AA^N$ induces an isomorphism from
$\grassS$ onto the closed subvariety of $\AA^N$ determined by the
following polynomials $\tau_q^\rho(X)$ in
$$K[X_{\alpha p ,q} \mid (\alpha,p)\ 
\sigma\text{-critical}, q \in \sigma(\alpha, p)].$$
Let $\frakA$ be the
polynomial ring with coefficients in $KQ$ (noncommutative in general) in
the variables $X_{\alpha p, q}$, where $(\alpha, p)$ ranges through the
$\sigma$-critical pairs, and $q$ through $\sigma(\alpha, p)$.  Moreover,
let $\frakC$ be the left ideal of $\frakA$ generated by the idempotents
$e_{\t+1}, \dots, e_n \in KQ$ annihilating $T$ together with the
differences
$$\alpha p - \sum_{q \in \sigma(\alpha,p)} X_{\alpha p, q} q.$$ 
Then $\frakA/\frakC$ is a free left module over the commutative
polynomial ring
$K[X_{\alpha p, q}]$ with basis $(q + \frakC)_{q \in \sigma}$.  In
particular, each element $\rho$ in $KQ$ is congruent, modulo $\frakC$, to
a sum $\sum_{q \in \sigma} \tau^{\rho}_q (X) q$ for unique polynomials
$\tau^{\rho}_q (X)$ in $K[X_{\alpha p, q}]$.  The affine vanishing set in
$\AA^N$ of the $\tau^{\rho}_q (X)$, where $\rho$ traces any left ideal
generating set $\R$ of $I$, depends only on $\S$, not on the
choice of $\R$.  As announced, this vanishing set is isomorphic to
$\grassS$ and will henceforth be identified with the latter variety
(see \cite{\class, Section 6} for a proof).   

We illustrate the procedure for finding the described
affine charts of $\grasstd$ with an example. Throughout,
$\seq$ denotes congruence modulo $\frakC$ in $\frakA$.  

\definition{Example 2.2} Let $\la = KQ /I$, where $Q$ is the quiver 

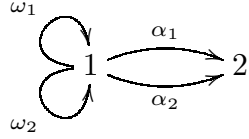
\ignore
$$\xymatrixrowsep{1pc}\xymatrixcolsep{3pc}
\xymatrix{ 1 \ar@(l,u)^{\omega_1} \ar@(l,d)_{\omega_2}
\ar@/^/[r]^{\alpha_1}
\ar@/_/[r]_{\alpha_2} &2  }$$
\endignore

\noindent and $I \subseteq KQ$ the ideal generated by
$\omega_i \omega_j$ for $i,j \in \{1,2\}$ and $\alpha_1 \omega_1 -
\alpha_2 \omega_2$.  Clearly, the listed relations generate $I$ as a left
ideal of $KQ$.

For $T = S_1$, we consider the 
$4$-dimensional abstract skeleton $\S = \{e_1, \omega_1,
\alpha_1\omega_1, \alpha_2\}$ with top $T$, compatible with the semisimple
sequence $\SS = (S_1, S_1 \oplus S_2, S_2)$.  This skeleton 
can be visualized as the set of paths that start at the top of the
following diagram:

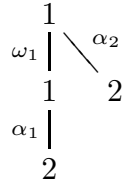
\ignore
$$\xymatrixrowsep{1.0pc}\xymatrixcolsep{0.8pc}
\vcenter{\xymatrix{ 1 \edge[d]_{\omega_1} \edge[dr]^{\alpha_2} \\ 1
\edge[d]_{\alpha_1} &2 \\ 2  }}$$
\endignore

\noindent In order to compute
$\grassS$, we list the
$\S$-critical pairs, together with their sets $\S(\alpha,p)$.  Obviously,
there is no harm in ignoring pairs
$(\alpha, p)$ with $\alpha p \in I$; the pairs that are left are
$(\omega_2, e_1)$ with $\S(\omega_2, e_1) = \{\omega_1\}$, 
$(\alpha_1, e_1)$ with $\S(\alpha_1, e_1)  = \{\alpha_2, \alpha_1
\omega_1\}$,  and $(\alpha_2, \omega_1)$ with
$\S(\alpha_2, \omega_1) = \{\alpha_1 \omega_1\}$.  This leads to the
following basic equivalences in $KQ[X_1, \dots, X_4]$, where
$X_1,\dots, X_4$ stand for $X_{\omega_2,\omega_1}, X_{\alpha_1,
\alpha_2}, X_{\alpha_1, \alpha_1 \omega_1}, X_{\alpha_2 \omega_1,
\alpha_1 \omega_1}$;  namely 
$\omega_2 \seq X_1 \omega_1$, $\alpha_1 \seq X_2 \alpha_2 + X_3
\alpha_1 \omega_1$, and
$\alpha_2 \omega_1 \seq X_4 \alpha_1
\omega_1$.   As one  finds (by an easy algorithm
described in \cite{\class, proof of Proposition 3.12}), all relations
$\omega_i \omega_j$ are congruent to zero modulo $\frakC$ in $\frakA$,
while
$\alpha_1 \omega_1 - \alpha_2
\omega_2$ is congruent to $(1 - X_4 X_1)\alpha_1\omega_1$.  Thus,
$$\grassS = \{(c_1,c_2,c_3,c_4) \in \AA^4 \mid 1 - c_1c_4 = 0\} \cong V(XY
- 1) \times \AA^2.  \qed$$   
\enddefinition

\head 3.  General facts about top-stable and layer-stable degenerations
\endhead

If $M$ and $M'$ are $d$-dimensional
$\la$-modules, represented by points $x$ and $x'$ in the classical affine
variety
$\modlad$, respectively, then $M'$ is said to be a {\it degeneration\/} of
$M$ if $x'$ belongs to the closure of the orbit $\GL_d.x$ in $\modlad$. 
Following \cite{\BongAdv}, we write $M
\degen M'$ in that case.  This notation reflects the fact that $\degen$ is
a partial order on the isomorphism classes of $d$-dimensional modules. 
Moreover, we say that $M'$ is a {\it top-stable degeneration\/} of $M$ if
$M \degen M'$ and $M/JM = M'/JM'$, and that $M'$ is a {\it layer-stable
degeneration\/} of $M$ if $M
\degen M'$ and $\SS(M) = \SS(M')$.  The latter means that
$J^lM/J^{l+1}M = J^lM'/J^{l+1}M'$ for all $l$.  (Recall that we identify
isomorphic semisimple modules.)

It has been the object of extensive efforts to characterize the
partial order
$\degen$ in purely algebraic terms.  This order is known to be
trapped between two algebraically defined partial orders which are
easier to track, namely the ``$\hom$-order", weaker than
$\degen$ in general, given by ``$M
\le_{\kleinhom} N \iff \dim_K \Hom_\la(M,X) \le \dim_K \Hom_\la(N,X)$ for
all
$X$" and the ``extension-order", stronger than $\degen$, arising as the
transitive closure of the relation ``$M \le_{\kleinext} N$ if $N \cong U
\oplus V$ for some exact sequence $0 \rightarrow U
\rightarrow M \rightarrow V\rightarrow 0$" (see, e.g., \cite{\Rie}, or
\cite{\Bongtrond} for an overview).  In case
$\la$ has finite representation type, the degeneration order has been
proved to coincide with the $\hom$-order; the final step in establishing
this is due to Zwara
\cite{\Zwa} and is the culmination of a line of
research going back to the 1980's; see e.g\.
\cite{\BongAdv, \Rie, \RieScho}.  In general, however, the two 
algebraically defined orders $\le_{\kleinhom}$ and
$\le_{\kleinext}$ yield only a loose bracket framing
$\degen$.  The following general description of $\degen$, again
due to Zwara \cite{\ZwaII} with roots in a result of Riedtmann
\cite{\Rie}, is in algebraic terms, but hard to track in general:
Namely, $M \degen M'$ precisely when there exists an exact sequence of
one of the following types:  
$0 \rightarrow M' \rightarrow M\oplus X \rightarrow X \rightarrow 0$ or
$0 \rightarrow X \rightarrow M\oplus X \rightarrow M' \rightarrow 0$ for
some module $X$.   Note that the implication ($M \le_{\kleinext} M'
\implies M
\degen M'$) specializes to $M \degen U \oplus M/U$ for every submodule
$U$ of $M$, and thus gives rise to the most accessible degenerations. 
This places a spotlight on indecomposable degenerations. 
\medskip

Fix $d\in\NN$, and a semisimple $\la$-module $T$ with projective cover $P$.
{\it None of the basic observations on layer- or top-stable degenerations
assembled in this section requires squarefreeness of
$T$.\/}  

Suppose $M$ and $M'$ are represented by points $x$
and $x'$ in $\toptd$.  Clearly,  
$M'$ is a top-stable degeneration of
$M$ precisely when $x'$ belongs to the relative closure of the orbit
$\GL_d.x$ in $\toptd$, and $M'$ is a layer-stable degeneration of $M$
precisely when $x'$ belongs to the closure of
$\GL_d.x$ in the subvariety $\layS$.  By the two-way
transfer of geometric properties between the
$\GL_d$-stable subsets of $\toptd$ on one hand and the
$\autlap$-stable subsets of $\grasstd$ on the other (cf\. Section 2),
this description of the top-, resp\. layer-stable, degenerations can be
reformulated as follows:  
 
\proclaim{Observation 3.1}  If $M \cong P/C$ with $C \in \grasstd$, then
a module $M'$ is a top-stable {\rm (}resp., layer-stable{\rm )}
degeneration of
$M$ if and only if $M' \cong P/C'$, where $C'$ belongs to the relative
closure of the orbit $\autlap.C$ in $\grasstd$ {\rm (}resp., in
$\grassSS${\rm)}.
\qed
\endproclaim

The next observation records some general facts concerning the
acting group $\autlap$ and its orbits.  They were established in
\cite{\class, Observation 2.3, Proposition 2.9, Lemma 2.10, and Proposition
4.1}; the structure of the $\unirad$-orbits arises as a special case of a
result of Rosenlicht \cite{\Ros} for arbitrary unipotent group actions. 
We start by noting  that  the unipotent radical of $\autlap$ is
$$\unirad = \{ \id + g \mid g \in \hom_\la(P,JP) \}.$$ 

\proclaim{Observation 3.2} Abbreviating $\unirad$ to $\U$, we obtain:  

$\bullet$  $\autlap \cong \U \rtimes \aut_\la(T)$ as algebraic groups.  In
particular,
$\autlap$ is a rational variety.

$\bullet$ For any point $C \in \grasstd$, say $C \in \grassSS$, the
$\U$-orbit $\U.C$ is closed in $\grassSS$ and isomorphic to the
affine space $\AA^{\m}$, where 
$$\m =\dim_K \hom_\la (P, JP/C) - \dim_K \hom_\la(P/C,JP/C).$$

$\bullet$  For any $C \in \grasstd$, the full orbit $\autlap.C$ has
dimension
$$\dim_K \hom_\la(P,P/C) - \dim_K \End_\la(P/C).$$

$\bullet$  For any $C \in \grasstd$, the orbit map $\autlap \rightarrow
\autlap.C$ is separable, and hence is a
geometric quotient of 
$\autlap$ by the stabilizer subgroup $\stab_{\autlap} C$.   
Restriction to $\U.C$ shows that the orbit map $\U \rightarrow
\U.C$ is a geometric quotient of $\U$ by $\stab_{\U} C$.  \qed
\endproclaim   

The following straightforward observation supplements the remark
about degenerations arising from short exact sequences.  Given
$M \in \lamod$, we call a submodule $U$ of $M$ {\it top-stably embedded\/}
in $M$ if $JM \cap U = JU$, {\it layer-stably embedded\/} in
$M$ in case $(J^l M) \cap U = J^l U$ for all $l \ge 1$.  Clearly, the
latter amounts to the same
as to say that the embedding
$U \hookrightarrow M$ induces monomorphisms $J^l U /J^{l+1} U
\hookrightarrow J^l M /J^{l+1}M$ for all $l \ge 0$, whence our
terminology. 

\proclaim{Observation 3.3}  Suppose $U$ is a submodule of $M$.  Then $U$ 
is top-stably {\rm (}resp., layer-stably{\rm )} embedded in
$M$ if and only if $U \oplus \bigl( M/U \bigr)$ is a top-stable {\rm
(}resp., layer-stable{\rm )} degeneration of $M$.  \qed
\endproclaim

In particular, the refined extension order, selecting short exact
sequences that remain exact on multiplication by powers of $J$, yields a
partial order on the isomorphism classes of layer-stable degenerations of
$M$.  

The corollary to the next observation provides an initial piece of
structural information on the top-stable degenerations of $M$ in comparison
with $M$: Roughly speaking, the degeneration process involves only upward
mobility of simple composition factors relative to the radical layering
of $M$.  The following partial order on the
$d$-dimensional semisimple sequences with top $T$ allows us to make this
precise. Namely, we define
$\SS \le
\SS'$ if and only if
$\bigoplus_{l \le r} \SS_l$ is a direct summand of $\bigoplus_{l \le r}
\SS'_l$ for all $r$.

The observation tying this order to degenerations is inspired by the
Schubert cells in the classical Grassmannian $\Gr(d', JP)$
containing $\grasstd$, in that its proof is based on intersection
dimensions with partial flags of subspaces of
$JP$.     

\proclaim{Observation 3.4}  Suppose that $\SS$ is a $d$-dimensional
semisimple sequence with top $T$.  Then the union $\bigcup_{\SS'
\ge \SS} \grass{\SS'}$ is closed in $\grasstd$. \endproclaim  

\demo{Proof} We write $J^l P / J^{l+1} P = \bigoplus_{1 \le i \le n}
S_i^{t_{li}}$ and $\SS_l =
\bigoplus_{1 \le i \le n} S_i^{s_{li}}$ for each $l \in \{1, \dots, L\}$; 
note that $s_{li} \le t_{li}$ by definition of semisimple sequences with
top $T$.  Moreover, for each vertex $e_i$ of
$Q$, we consider the following partial flag of subspaces of the $K$-space
$JP$:
$$e_i J^L P \subseteq e_i J^{L - 1} P \subseteq \cdots
\subseteq e_i J P.$$ 
Then the union of the $\grass{\SS'}$, with $\SS'$
tracing the semisimple sequences larger than or equal to
$\SS$ under $\le$, coincides with $\grasstd \cap 
\bigl(\bigcap_{1 \le i \le n} V_i\bigr)$, where $V_i$ is the following
subset of the classical Grassmannian $\Gr(d',JP)$:
$$V_i = \{C \in \Gr(d',JP) \mid \dim \bigl( C \cap e_i J^l P \bigr) 
\ge  \sum_{k \ge l}
\bigl( t_{ki} - s_{ki} \bigr) \text{\ for all\ } 1 \le l \le L\}.$$ 
The
latter sets are well known to be closed in $\Gr(d',JP)$, and our claim
follows. \qed \enddemo 

\proclaim{Corollary 3.5}  Suppose $M'$ is a top-stable degeneration of
$M$.  Then $\SS(M') \ge \SS(M)$.  If $M'$ fails to be a
layer-stable degeneration of $M$, this inequality is strict. \qed
\endproclaim

Finally, we note that we can always obtain the top-stable
degenerations $P/C'$ of a module $P/C$ along projective curves, due
to completeness of the orbit closure of $C$.  Upgrades of this
observation for simple, resp\. squarefree, top $T$ will be given in
Proposition 4.2 and Corollary 5.4, respectively.  The usefulness of this
approach, both towards proofs and concrete examples, lies in the fact that
the geometry of the
$\autlap$-orbits of $\grasstd$ is particularly transparent in these
cases.  Moving subspace flags of the point $C \in \grasstd$ along
in tandem, will moreover enable us to link up the structures of
$P/C$ and $P/C'$.       

Let $\flagtd \subseteq \Gr(1,JP) \times \Gr(2,JP) \times \cdots
\times \Gr(d',JP)$ denote the variety of subspace flags of the
$d'$-dimensional $\la$-sub{\it modules\/} of $JP$, where $d' = \dim P -
d$.  In other words, if $\F$ is the full flag variety of the $K$-space
$JP$, and $\pi$ denotes the
projection from $\Gr(1,JP) \times \cdots \times \Gr(\dim JP, JP)$ onto 
$\Gr(1,JP) \times  \cdots \times \Gr(d',JP)$, then $\flagtd$ equals the
intersection
$$\pi(\F) \ \cap \ \biggl(\Gr(1,JP) \times \cdots \times \Gr(d' - 1,
JP) \times \grasstd\biggr).$$
We observe that this intersection is closed in $\Gr(1,JP) \times  \cdots 
\times \Gr(d',JP)$, since $\grasstd$ is closed in $\Gr(d',JP)$ and the
variety $\F$  --  whence also $\pi(\F)$  --  is complete.  Thus $\flagtd$,
in turn, is projective.   Moreover, we note that the action of $\autlap$
on $\grasstd$ naturally extends to an $\autlap$-action on
$\flagtd$.

The following result is due to J\'anos Koll\'ar (personal communication).  It
extends to higher generality than required for our application; here we are
only concerned with classical projective varieties, namely vanishing sets of
families of homogeneous polynomials in projective space.  Recall that an
irreducible variety
$V$ is called unirational if its function field embeds into a finitely
generated purely transcendental extension of the base field $K$, the latter
amounting to the existence of a dominant rational map from some affine space
$\AA^r$ to $V$ (see, e.g., \cite{\Bor, Chapter AG, 13.7}).  The analogous
result for characteristic zero (based on Hironaka's resolution of
singularities) was already obtained and used by Hanspeter Kraft in
the context of degenerations of finite dimensional representations; see
\cite{\Kraft, p\. 227}.

\proclaim{Proposition 3.6} {\rm{[J. Koll\'ar]}} Let $V$ be a
unirational projective variety of positive dimension {\rm{(}}over an
algebraically closed field $K$ of arbitrary characteristic{\rm{)}} and
$c_1, \dots, c_s \in V$ any finite collection of points.  Then there exists a
curve $\psi: \PP^1 \rightarrow V$, the image of which contains all $c_i$. 
\endproclaim

\demo{Proof}  We fix an embedding $V \subseteq \PP^R$ of $V$ into a
projective space.  By unirationality of $V$, there exists a dominant rational
map from some
$\PP^r$ to $V$, represented by $(F,U)$ say, where $U$ is a dense open
subset of $\PP^r$ and 
$$F(x) = \bigl(f_0(x), \dots, , f_R(x) \bigr)$$ 
with homogeneous polynomials $f_j \in K[X_0, \dots, X_r]$.  First we cut the
variety
$V$ with general hypersurfaces of $\PP^R$ containing $c_1, \dots, c_s$,
to obtain an irreducible curve $\C$ in $V$ with $\{c_1, \dots, c_s\} \subset
\C$ such that $F$ induces a birational equivalence between $\C$ and a
suitable projective curve $\C'$ in $\PP^r$.  (Feasibility is a consequence
of a Bertini-type theorem, which is known to hold in arbitrary
characteristic:  Namely, if
$X$ is an irreducible subvariety of dimension at least $2$ of a projective
space $\PP^R$, then, generically, hypersurface sections $X \cap H$ in $\PP^R$
are again irreducible; see \cite{\Jou}.  In fact, the result has been
generalized to not necessarily complete linear systems of divisors without
fixed components on
$X$, so as to yield irreducibility of their generic
members.  In our specific situation, one applies this version of Bertini's
Theorem to the linear system of hyperplanes of a fixed degree
$\ge s+1$ containing the points $c_1, \dots, c_s$ of $V$.)   Let
$h:  \D \rightarrow \C'$ be the normalization of $\C'$.  This provides us
with a smooth projective curve $\D$ with the property that 
$$\{c_1, \dots, c_s\} \subset \overline{(F\circ h) \bigl(h^{-1}(U\cap \C')
\bigr)} \subseteq V.$$
Since $\D \setminus h^{-1}(U \cap \C')$ is finite, smoothness allows us to
extend $F \circ h$ to a morphism $H: \D \rightarrow V$ (see, e.g.,
\cite{\Har, Chapter I, Proposition 6.8}), which yields
$\{c_1, \dots, c_s\} \subset H(\D)$.  

Let $d_i \in \D$ with $H(d_i) = c_i$.  Our claim will follow from the
following two facts:

(1)  For every $m \ge 1$, there exists a morphism $g: \PP^1 \rightarrow
\PP^r$, together with points $b_1, \dots, b_s \in \PP^1$, such that the
Taylor expansion of $g$ at $b_i$ agrees, to order $m$, with the Taylor
expansion of  $h$ at $d_i$, for all $i$.

(2)  There is an integer $m \ge 1$ with the following property: If $g: B
\rightarrow \PP^r$ is any morphism from a smooth projective curve $B$ and if 
$b_1, \dots, b_s \in B$ are points such that the Taylor expansion of $g$
at
$b_i$ agrees, to order $m$, with the Taylor expansion of $h$ at $d_i$, for
all $i$, then $\{c_1, \dots, c_s\} \subset \overline{(F \circ
g)\bigl(B \cap g^{-1}(U)\bigr)}$.

For the first of these facts, we refer to \cite{\Kol, (4.1.2.4)}, where (1)
is stated for the base field $K = \CC$.  The proof given in
\cite{\Kol, (5.2)} works in arbitrary characteristic, however.

The second can be seen as follows:  Suppose that $g : B \rightarrow
\PP^r$ is a morphism from a smooth curve $B$, and let $b_1, \dots, b_s \in
B$.  If $G : B
\rightarrow V$ is the extension of $F \circ \bigl(g|_{g^{-1}(U)}\bigr)$ to
$B$, then evidently $G(B) = \overline{(F \circ
g)\bigl(B \cap g^{-1}(U)\bigr)}$.  We focus on a point $d_i$, pick a local
coordinate $t_i$ for $\D$ at $d_i$, and consider $h$ as given by power 
series
$h_j$ in $t_i$ near $d_i$, that is, as given in the form $\bigl(h_0(t_i),
\dots, h_r(t_i)\bigr)$.  Then
$$H(t_i) = \bigl(f_0\bigl(h_0(t_i), \dots, h_r(t_i)\bigr), \dots,
f_R\bigl(h_0(t_i), \dots, h_r(t_i)\bigr) \bigr)$$
in a neighborhood of $d_i$, punctured at $d_i$ (the composition $F
\circ h$ need not be defined at $d_i$), alias, for $t_i$ near $0$ but
different from $0$.  Let
$t_i^{m_i}$ be the minimal order of vanishing of the
$f_j\bigl(h_0(t_i), \dots, h_r(t_i)\bigr)$.  Since 
$H(t_i) = (1/t_i^{m_i}) H(t_i)$ in $\PP^r$ whenever $t_i \ne 0$  
in said neighborhood, $H(0)$ is obtained by substituting $t_i = 0$ in the
right-hand side of the previous equality.
Thus, if $g_j(t_i)$ are power series expansions of the coordinate functions
of $g$  near
$b_i$ (where $t_i$ also stands for a uniformizing parameter at $b_i$ of the
smooth curve $B$) such that
$$f_j\bigl(h_0(t_i), \dots, h_r(t_i)\bigr) \equiv f_j\bigl(g_0(t_i), \dots,
g_r(t_i)\bigr) \quad \text{mod}\ t_i^{m_i+1},$$
for $1 \le j \le R$, then $H$ coincides with $G$ at $t_i = 0$, and
returning to the original variables, we obtain
$c_i = H(d_i) = G(b_i)$.   Hence $m \ge \max\{m_1, \dots, m_s\} + 1$
satisfies the requirement under (2).
\qed 
\enddemo

The first of the observations below parallels Kraft's use of curves, over the
base field $\CC$, to reach a degeneration $M'$ of a representation $M$ of a
complex algebra, by way of a rational curve connecting a point in the orbit
of $M$ in $\modlad$ to a point representing $M'$; he based his
result on the special case of Proposition 3.6 for characteristic
zero mentioned above.  In the framework of $\grasstd$, we have
the added advantage of a manageable format for such curves, improving the
theoretical and computational accessibility of degenerations.

\proclaim{Observation 3.7} {\rm (1)} For $C \in \grasstd$ and any
top-stable degeneration $P/C'$ of
$P/C$, there exists a curve $\psi: \PP^1 \rightarrow \overline{\autlap.C}$
such that $\psi^{-1}(\autlap.C)$ is dense in $\PP^1$ and $C' \in
\Im(\psi)$.

{\rm (2)}  If $\C$ is a subspace flag of $C$ and $\C'$ a
point in the closure of $\autlap.\C$ in $\flagtd$, there exists a curve 
$\psi: \PP^1 \rightarrow \overline{\autlap.\C}$
such that $\psi^{-1}(\autlap.\C)$ is dense in $\PP^1$ and $\C' \in
\Im(\psi)$.  

{\rm (3)} Suppose that $\C =
\bigl(C_1, \dots, C_{d'}\bigr) \in \flagtd$ is a flag with $C_{d'} = C$.
Moreover, let $\rho: U \rightarrow \autlap$ be a morphism defined on a
dense subvariety $U$ of
$\PP^1$.  Then the morphism $\psi: U \rightarrow \autlap.\C \subseteq
\flagtd$, defined by $\tau
\mapsto \bigl( \rho(\tau).C_1,\dots,\rho(\tau).C_{d'} \bigr)$, extends to a
unique morphism
$\overline{\psi}: \PP^1 \rightarrow \overline{\autlap.\C}$.  The latter
maps any point $\tau \in \PP^1$ to a flag $\D =
\bigl(D_1, \dots, D_{d'}\bigr)$ of the point $D = D_{d'}$ in the orbit
closure $\overline{\autlap.C}$.  In case $\C$ is a flag consisting
of $\la$-submodules of $C$, the flag
$\D$ consists of $\la$-submodules of $D$.
\endproclaim

\demo{Proof} (1) $\autlap$ being a rational variety by
Observation 3.2, each orbit closure is unirational.  Thus Proposition 3.6
guarantees that any two points of an orbit can be connected by a rational
curve.  In other words, there exists a morphism $\psi: \PP^1
\rightarrow \overline{\autlap.C}$,
such that $C, C' \in \Im(\psi)$.  Because $\autlap.C$ is locally closed in
$\grasstd$, our first claim follows.  The proof of (2) is analogous.

(3) In light of the fact that $\flagtd$ is projective, the extension
statement for curves of flags is immediate.   The final assertion under
(3) follows from the fact that 
$$\flagtd \cap \bigl(\frak{Grass}^T_{\dim P
-1} \times \cdots \times \frak{Grass}^T_{\dim P - d'}\bigr)$$ 
is a closed subvariety of $\flagtd$.  \qed
\enddemo

The first part of Observation 3.7 has an analogue for the
classical module variety:  Namely, given any point
$x \in \modlad$,  each element in the closure of
$\GL_d.x$ can be realized in the form $\psi(a)$ for some morphism
$\psi: U \rightarrow \overline{\GL_d.x}$, where $U$ is a nonempty open
subset of
$\AA^1$ containing $a$, and $\psi\bigl( U\setminus \{a\} \bigr)
\subseteq \GL_d.x$; see \cite{\Kraft} and \cite{\ZwaII}. 
The classical variety only sporadically lends itself to explicit
computations of all degenerations of a given module, however.  This is
otherwise in the projective setting, as evidenced in Sections 4.B and 5.B.

\definition{Definitions and Remarks 3.8} Following standard terminology,
we refer to a {\it proper\/} top-stable (resp., layer-stable) degeneration
$M'$ of
$M$ as {\it minimal top-stable\/}  (resp., {\it minimal layer-stable\/})
in case there is no top-stable (resp., layer-stable) degeneration lying
properly between $M$ and $M'$ in the degeneration order.  Then $M'$ is
minimal also in the poset of {\it all\/} proper degenerations of $M$,
since any degeneration of $M$ which is trapped between $M$ and $M'$ is in
turn top-stable (layer-stable).   

{\it Maximal\/} top- or layer-stable degenerations are defined analogously,
on waiving the condition that they be proper.  Clearly
a module
$M$ has no proper top-stable degenerations precisely when the set of
(isomorphism classes of) maximal top-stable degenerations consists of $M$
alone.  Whereas the unique maximal degeneration of
$M$, namely the direct sum of the simple composition factors of $M$, does
not hold much interest, the maximal {\it top-stable\/} and the maximal
{\it layer-stable\/} degenerations of
$M$ do;  in particular, we will see that, in general, even modules with
simple tops may have infinite families of maximal
top-stable degenerations requiring parameter spaces of arbitrarily high
dimension.

Finally, we say that $M'$ is a top-stable (layer-stable) degeneration of
$M$ {\it of height\/}
$h$, if
$h$ is the maximal length of a chain of top-stable (layer-stable)
degenerations connecting
$M$ to $M'$.  (As in the theory of prime ideals, we count as the Romans
did, calling $M = M_0 <_{\text{deg}} M_1 <_{\text{deg}} \dots
<_{\text{deg}} M_l$ a chain of length $l$;  thus the minimal top-stable
degenerations of
$M$ are those of height $1$.)  Clearly, the heights of the
top-stable (layer-stable) degenerations of a $d$-dimensional module $M =
P/C$ are bounded from above by the dimension of the orbit $\autlap.C$ in
$\grasstd$.
\enddefinition

\head 4. The local case: $T$ simple \endhead

Throughout this section, we assume $T$ to be simple, i.e., $T =
\la e/ Je$ for some primitive idempotent $e$, and $P = \la e$.  This means
that every module with top $T$ is local, i.e., has a unique maximal
submodule.  The local case affords several bonuses, sufficiently incisive
to warrant separate treatment, all the more since all maximal top-stable
degenerations of an arbitrary finite dimensional module are direct sums of
local ones. 

Responsible for these bonuses is the fact that, in this situation,
$\aut_\la(P)$ is the direct product
$K^* \times \unirad$, where the second factor is the unipotent radical of
$\autlap$ identified in Section 3, and the first acts trivially on
$\grasstd$.  In particular, the
$\autlap$-orbits in
$\grasstd$ coincide with the $\unirad$-orbits, and hence are affine
spaces, by Observation 3.2.  
\medskip
  
We start with a few elementary comments, the first giving an explicit
incarnation of an isomorphism between the varieties $\AA^\m$ and
$\autlap.C$, for $C\in \grasstd$.  Since $T = \la e / Je$, the group
$\unirad$ is anti-isomorphic to the multiplicative subgroup
$e + eJe$ of the group of units of $e\la e$,
via the assignment which sends any  $e+a$ to right multiplication of $P$ by
this element.   

Now suppose that $\omega_1, \cdots,
\omega_\mu$ form a basis for $eJe$.  We
choose the $\omega_i$ to be oriented cycles from $e$ to
$e$, which is clearly always possible.  Then the above isomorphism of
groups induces an isomorphism
$\bigoplus_{1 \le i \le \mu} K\omega_i \cong \unirad$ of algebraic
varieties, which maps any element $a = \sum_{1 \le i \le \mu} k_i
\omega_i$ in
$eJe$ to right multiplication on $P$ by $e + a$.  To describe the orbits
$\autlap.C$ in these terms, let
$\stab_{eJe} C$ be the $K$-space consisting of those elements
$a \in eJe$ for which $Ca \subseteq C$, and denote the unipotent radical
$\unirad$ by $\U$.  As a variety, $\stab_{eJe} C$ then becomes isomorphic
to $\stab_{\U} C$ under the described assignment, whence $\dim
\autlap.C$ equals $\dim_K \bigl( eJe/\stab_{eJe} C \bigr)$; this
dimension also equals the invariant $\m$ in the local case (again see
Observation 3.2).  Provided that 
$\omega_1,
\dots,
\omega_\m$ form a
$K$-basis of $eJe$ modulo
$\stab_{eJe} C$, this latter map in turn gives rise to an isomorphism
$\bigoplus_{1 \le i \le \m} K\omega_i \cong \U/ \stab_\U C$.  We thus
obtain an explicit isomorphism between $\autlap.C$ and $\AA^{\m}$ as
follows:

\proclaim{Lemma 4.1}  Let $C$ and $\omega_1, \dots, \omega_\m$ be as
above.  Then the map 
$$\bigoplus_{1 \le i \le \m} K\omega_i \rightarrow \autlap.C,$$ 
sending any
element $a$ in $\bigoplus_{1 \le i \le \m} K\omega_i$ to the submodule
$C(e + a)$ of $P$, is an isomorphism of varieties. \qed
\endproclaim

Based on the previous lemma, we can strengthen Observation 3.7.  In light
of projectivity of $\grasstd$, every morphism $\psi: U
\rightarrow \grasstd$, where $U$ is a dense subset of $\AA^1$, can be
uniquely extended to $U \cup \{\infty\} \subseteq \PP^1$.  Therefore the
notation $\lim_{\tau
\rightarrow \infty} \psi(\tau)$ is unambiguous. Moreover, closedness of
the flag variety $\flagtd$ in $\Gr(1,JP) \times \Gr(2,JP) \times
\cdots \times \Gr(d',JP)$ allows for an analogous limit notation for curves
in $\flagtd$; as before, $d' = \dim P - d$.

For the first assertion of Proposition 4.2, we only need to note that any
curve $U \rightarrow \autlap.C$ with $U \subseteq \AA^1$ dense
coincides, on a dense subset of $U$, with one of the form $\tau
\mapsto C \bigl(e + \sum_{1 \le i \le \m} r_i(\tau) \omega_i \bigr)$ for
suitable rational functions $r_i(\tau)$ (Lemma 4.1).  The final assertion
rests on the fact that all $\autlap$-orbits of $\grasstd$ are affine.  

\proclaim{Proposition 4.2}  {\rm (1)}  Let $C, C' \in \grasstd$ be such 
that
$P/C$ degenerates to $P/C'$, and let $\omega_1, \dots, \omega_m$
denote oriented cycles $e \rightarrow e$ in $Q$ generating the $K$-space
$eJe$ modulo
$\stab_{eJe} C$.  Then there exist polynomials $q(\tau)$ and $p_i(\tau)$ 
in $K[\tau]$ for $i \le m$ with $q \ne 0$ such that 
$$\lim_{\tau \rightarrow \infty} C
e(\tau) = C',$$
 where, for each $\tau \in K$,  
$$e(\tau) = q(\tau) e + \sum_{1 \le i \le m} p_i(\tau) \omega_i$$ 
is
identified with an element of $e\la e$, a unit in $e \la e$ whenever
$q(\tau) \ne 0$.

Moreover, given any flag $\C = (C_j)_{1 \le j \le d'}$ of subspaces $C_j$
of $C$, the ``limit" 
$$\C' = \lim_{\tau \rightarrow \infty} \C e(\tau) = \bigl(\lim_{\tau
\rightarrow \infty} C_j e(\tau) \bigr)_{j \le d'}$$  is a flag of $C'$. 
In particular, any subspace $D$ of $C$ gives rise to a
{\rm (}$\dim D${\rm )}-dimensional subspace of $D' = \lim_{\tau \rightarrow
\infty} D e(\tau)$ of $C'$.  If $D$ is a $\la$-submodule of $P$, then so is
$D'$. 
\smallskip

{\rm (2)} Conversely, given arbitrary polynomials $q(\tau)$ and
$p_i(\tau)$ with $q \ne 0$, we set
$$e(\tau) =   q(\tau) e
+ \sum_{1 \le i \le m} p_i(\tau) \omega_i \quad \text{for all\ } \tau
\in \AA^1.$$  
If the curve $\tau \mapsto C e(\tau)$ in
$\grasstd$, defined on the complement of the zero set of $q$ in $\AA^1$, 
is nonconstant, its unique extension to
$\PP^1$ takes a value
$C'$ such that $P/C'$ is a proper degeneration of $P/C$. \qed
\endproclaim
 
A proof of the following lemma has been communicated to the author by Ron
Donagi; we include it with his permission. 

\proclaim{Lemma 4.3} {\rm{[R\. Donagi]}}  Let $V$ be an $m$-dimensional
irreducible projective variety and $W$ a nontrivial open subvariety
permitting a non-constant regular function $W \rightarrow \AA^1$.  Then $V
\setminus W$ has dimension
$m-1$. 

In particular, $\dim V \setminus W = m-1$ whenever $V$ is a projective
variety containing a dense open subset $W$ isomorphic to $\AA^m$.
\endproclaim

\demo{Proof} That $\dim\bigl(V \setminus W\bigr) \le m-1$ is clear.  

For
the converse, let  $\phi:  W \rightarrow \AA^1$ be a nonconstant regular
function as postulated, and let $\Gamma \subseteq V \times \PP^1$ be the
closure of its graph in $V \times \PP^1$.  Moreover, let
$\overline{p}$ and
$\overline{q}$ be the projections from $V \times \PP^1$ onto $V$ and
$\PP^1$, respectively, and $p$, $q$ their restrictions to $\Gamma$.  Due
to projectivity of $\Gamma$ and our assumption on $\phi$, the image of $q$ is
$\PP^1$, and therefore the fiber $q^{-1}(\infty) \subseteq \Gamma$ has
dimension at least $\dim \Gamma - 1 = m - 1$ by \cite{\Bor, Chapter AG,
Theorem 10\.1}).  Clearly, the restriction of $p$ to $q^{-1}(\infty)$ is
injective; indeed this restriction is just the embedding of
$q^{-1}(\infty)$ into $\Gamma \cap \bigl(V\times \{\infty\}\bigr)$ followed
by the isomorphism $\overline{p}: V \times \{\infty\} \rightarrow V$. 
Consequently, the (closed) subvariety $p\bigl(q^{-1}(\infty)\bigr)$ of $V$
has dimension at least
$m-1$ as well.  Clearly, it is contained in $V \setminus W$, and the
missing inequality follows. \qed    
\enddemo

\subhead 4.A. The basic theorems for the local case \endsubhead

We have seen that all $\autlap$-orbits of $\grasstd$ are full affine
spaces in the local case.  For the sake of symmetry  --  compare 
with the more general statement addressing the geometry of the
$\autlap$-orbits in the twin Theorem 5.1  --  we include this
fact as part (1) of the theorem below.  For a proof of our first
theorem we only need to assemble the information accumulated so far. In
the subsequent theorem, we provide more detail on the possible sizes and
shapes of the posets of top-stable degenerations. This second result,
proved in Section 4.B, demonstrates that our conjectures concerning the
general form of the poset of top-stable degenerations of a local module
cannot be simplified.   

\proclaim{Theorem 4.4}   Suppose $T = \la e / Je$ is simple and $P = \la
e$.
\medskip

\noindent {\rm \bf{(1)} Structure of the $\autlap$-orbits and chain
lengths.} For each $C\in \grasstd$, 
$$\autlap.C \cong \AA^{\m} \qquad \text{where} \ \ \m = \dim_K
\hom_\la (P, P/C) - \dim_K \End_{\la}(P/C).$$
In particular, the lengths of chains
of top-stable degenerations of $P/C$ are bounded above by $\m$. 
\medskip

\noindent {\rm \bf{(2)} Layer-stable degenerations.} Local modules do not
have any proper layer-stable degenerations.
\medskip

\noindent {\rm \bf{(3)} Top-stable degenerations.} Let $M$ be a local
module with top $T$, say $M \cong P/C$ with $C \in \grasstd$.
\smallskip

{\rm \bf{(a)} Existence.}  The following conditions are equivalent:
\roster
\item"(i)" $M$ does not have any proper top-stable degenerations.

\item"(ii)" The orbit $\autlap.C$ is a singleton.

\item"(iii)" $\dim_K \End_{\la}(M)$ equals the multiplicity of $T$ as a
composition factor of $M$.

\item"(iv)"  $C$ is invariant under all endomorphisms of $P$, that is,
$C\omega \subseteq C$ for all oriented cycles
$\omega$ from $e$ to $e$ in the quiver $Q$.
\endroster

\noindent In particular, these equivalent conditions are satisfied if all
appearances of $T$ as a composition factor of $JM$ occur in the  socle
of $JM$.  
\smallskip

{\rm \bf{(b)} Unique existence.}  
$M$ has a unique proper top-stable degeneration if and only if $\m = 1$,
that is, if and only if
 the multiplicity of the simple module $T$ in $M$ exceeds $\dim_K
\End_\la(M)$ by $1$.
\smallskip

{\rm \bf{(c)} Bases.}  If $M'$ is a top-stable degeneration of $M$, then
$M$ and $M'$ share a basis consisting of paths.  More precisely, there
exist submodules $D$ and $D'$ of $P$ with $M \cong P/D$ and $M' \cong
P/D'$, together with a set
$\B$ of paths in $KQ$ {\rm{(}}which is closed under right subpaths{\rm{)}},
such that
$\{ p + D \mid p \in \B \}$ is a basis for $P/D$ and $\{ p + D' \mid p \in
\B \}$ a basis for $P/D'$.  
\smallskip

{\rm \bf{(d)} The maximal top-stable degenerations of $M$} always possess
a fine moduli space classifying them up to isomorphism, namely 
$$\maxtopdeg(M)  = \{C' \in \overline{\autlap.C} \mid \autlap.C' \text{\
is a singleton}\},$$ and the latter is a projective variety of dimension
at most $\max\{0, \m - 1\}$.  More specifically, the assignment
$C' \mapsto P/C'$ yields a universal {\rm(}bijective{\rm)} parametrization
of the maximal top-stable degenerations of $M$, up to isomorphism.  
\endproclaim

As a special case of part (2), we rediscover the fact that uniserial
modules fail to have proper uniserial degenerations (see \cite{\Bong}).
Moreover, we note that the conditions in Theorem 4.4 characterizing
existence and unique existence of proper top-stable degenerations can be
readily checked in concrete situations, given $C$ and a quiver
presentation of $\la$; see Section 4.B.  

\demo{Proof} For (1), see Lemma 4.1, and for (2), see \cite{\class,
Corollary 4.3}.

(3a) The equivalences can be readily derived from
\cite{\class, Theorem 4.2}; we include the short argument for the
convenience of the reader.  In view of part (1), the equivalence of
statements (ii)--(iv) is straightforward, as is the implication `(ii)
$\implies$ (i)'.  For the converse, we only need to combine Observation
3.1 with the fact that the closure of $\autlap.C$ in $\grasstd$ is
projective.  The supplementary statement under (3a) is an obvious
consequence of the equivalences.

(3b) First suppose that $\autlap.C$ is one-dimensional.  That this
ensures  existence of a proper top-stable degeneration of
$P/C$ follows from part (a). Moreover, since
$\autlap.C \cong \AA^1$ by part (1), the closure of
$\autlap.C$ in the projective variety $\grasstd$ contains only a single
point outside $\autlap.C$:  Indeed, let $\AA^1 \cong \autlap.C \rightarrow
\overline{\autlap.C}$ be the embedding. By \cite{\Har, Chapter I,
Proposition 6.8}, this embedding can be extended to a morphism $\PP^1
\rightarrow \overline{\autlap.C}$, and by \cite{\Bor, Chapter AG, Section
7.4}, the image of $\PP^1$ under this extension is closed, whence
$\overline{\autlap.C} \setminus
\autlap.C$ is a singleton.  This guarantees uniqueness.  The converse is due
to the  fact that, for
$\dim \autlap.C \ge 2$, the variety
$\overline{\autlap.C} \setminus \autlap.C$ has dimension at least $1$;
combine Lemma 4.3 with part (1) for verification.  Thus, if this set
difference does not contain an orbit of positive dimension, it contains 
at least two distinct orbits.  If, on the other hand, it does contain an
orbit of positive dimension,
$\autlap.C'$ say, then
$P/C'$ in turn has a proper top-stable degeneration by part (a), which
again yields two distinct orbits in the boundary of $\autlap.C$.      

(3c) Let $M'
\cong P/D'$, and pick any skeleton $\sigma$ of $M'$.  Then $D'$ belongs to
$\grassS$ and hence to the open subset 
$$\aux(\sigma) = \{E \in \grasstd \mid P/E \text{\ has basis\ } \{p + E
\mid p \in \sigma\} \}$$  
of $\grasstd$.  Since $M \cong P/C$, this means that the intersection
$\aux(\sigma) \cap
\overline{\autlap.C}$ is nonempty and open in $\overline{\autlap.C}$. 
The orbit $\autlap.C$ being open in its closure and irreducible by part (1),
we conclude that $\aux(\sigma) \cap \autlap.C \ne \varnothing$, which yields
a point $D \in
\autlap.C$ such that $P/D$ has basis $\{p + D \mid p \in \sigma\}$ as
required; see Section 2.

(3d) That $\maxtopdeg(M)$ is a fine moduli space for the maximal
top-stable degenerations of $M$ follows from (3a) and \cite{\class, Theorem
4.4}. Projectivity of $\maxtopdeg(M)$ follows from closedness in
$\overline{\autlap.C}$, which, in turn, is due to the lower semicontinuity
of the orbit dimension.  Since $\dim \autlap.C = \m$ and $\maxtopdeg(M)$
is either a singleton or contained in $\overline{\autlap.C} \setminus
\autlap.C$, the bound on the dimension is obvious.
\qed  
\enddemo

Next, we provide evidence for the fact that the size of the poset of
top-stable degenerations of a local module $M = P/C$ may be very large,
both ``horizontally" and ``vertically", within the limits imposed by the
equality $\dim \overline{\autlap.C} \setminus \autlap.C = \dim
\autlap.C - 1$ (Lemma 4.3).

\proclaim{Theorem 4.5}  Let $n_1, \dots, n_s$ be arbitrary 
positive integers, and $l \ge s$. Then there is a finite dimensional
algebra
$\la$, together with a local $\la$-module $M = P/C$, such that the poset
of top-stable degenerations of $M$ has the following properties: 
$$\maxtopdeg(M) \cong
\PP^{n_1} \times \dots \times \PP^{n_s},$$   
the integer $l$ is the minimal length of a saturated chain of top-stable
degenerations of $M$, and $\dim \autlap.C = l + \sum_{1 \le i \le s}
n_i$.   

In addition, we may require that either

{\rm (i)}  the lengths of the saturated chains of top-stable degenerations
of $P/C$ trace all values between $l$ and $(l-s)+ \sum_{1 \le i \le s}
n_i$,
or else

{\rm (ii)}  the saturated chains of top-stable degenerations of $P/C$
have constant length $l$.          
\endproclaim

Theorem 4.5 will be established in Section 4.B, Examples 4.10.  Beyond the
listed properties, our families of examples exhibit the following feature,
which we believe to be representative of the general situation, in
a sense made more precise below:  For every value $h \le \dim \autlap.C$
occurring as the height of a top-stable degeneration of $P/C$, the set of
(isomorphism classes of) top-stable degenerations of that height is the
disjoint union of finitely many families, each of which is bijectively
parametrized by some torus orbit in a direct product of projective spaces. 
By a {\it torus orbit\/} in a product $\PP^{n_1} \times \cdots \times
\PP^{n_s}$ we will mean any orbit of the canonical
$(K^*)^{n_1+1} \times
\cdots \times (K^*)^{n_s+1}$-action.

We conclude the section with a number of conjectures, setting up as tight
a framework for classification of the poset of top-stable degenerations of
local modules as can be expected in this generality.  The
first conjecture pertains to modules with arbitrary squarefree tops, the
second to local modules.  

\definition{Conjectures 4.6} 
 For any module $M = P/C$ with squarefree top, the following hold:
\smallskip

{\bf{Conjecture (1)}} There exist finitely many non-negative integers
$n_1 \le \cdots \le n_s$, repetitions allowed, such that the
following conditions are satisfied:

(a)  The irreducible components of the fine moduli space $\maxtopdeg(M)$
coincide with the connected components, and each component is
isomorphic to a direct product of projective spaces, $\prod_{i
\in I} \PP^{n_i}$, for a suitable subset $I$ of $\{1, \dots, s\}$.  

(b) For each $h \le \dim \autlap.C$, the top-stable degenerations of height
$h$ above $M$ can be partitioned into finitely many disjoint families, each
of which is bijectively parametrized by a torus orbit in a product 
$\prod_{i \in I} \PP^{n_i}$ for some $I \subseteq \{1, \dots, s\}$.
\smallskip

If $s$ and $\sum_{1 \le i \le s} n_i$ are minimal with respect to (a) and
(b), then $n_1, \dots, n_s$  are isomorphism invariants of $M$, and 

(c) the difference 
$\dim \autlap.C - \sum_{1 \le i \le s} n_i$ coincides with the minimum of
the lengths of saturated chains of top-stable degenerations of $M$.
\smallskip

{\bf{Conjecture (2)}} Now suppose that $M$ has simple top, and let $\m =
\dim \hom_\la(P,JM) - \dim \hom_\la(M,JM)$.  Then
$\maxtopdeg(M)$ is irreducible, that is, the universal family of
isomorphism classes of maximal top-stable degenerations of
$M$ is a $\bigl(\prod_{i \in I} \PP^{n_i}\bigr)$-family for a suitable
subset $I \subseteq \{1, \dots, s\}$.  If, in addition, the poset of
(isomorphism types of) top-stable degenerations of
$M$ is finite, then its cardinality (counting $M$) is bounded from below
by $\m+1$ and from above by $2^{\m}$.  More precisely, for $0 \le j \le
\m$, the number of top-stable degenerations of $M$ of height $j$ lies
between
$1$ and $\m \choose j$. In particular, this means that $M$ has a unique
maximal top-stable degeneration, and that the latter has height $\m$. 
Moreover, all saturated chains linking
$M$ to its unique maximal top-stable degeneration have the same
length.  (This catenarity condition fails when
either of the hypotheses is removed.)   
\enddefinition

Suppose $M$ is a module with squarefree top and $s, n_1, \dots, n_s$
are as postulated in  Conjecture (1) above, chosen to be  minimal.  The
requirement that
$M$ have only finitely many top-stable degenerations then becomes
equivalent to
$s = 1$ and $n_1 = 0$; for an infinite poset of top-stable degenerations,
the integers
$n_i$ are all positive, due to minimality.  Provided that the conjectures
can be confirmed, all bounds are sharp.  This will be seen in Section 4.B
as an application of the theory developed so far.  In all  of the local
examples given here,
$\dim \maxtopdeg(M) = n_1 + \cdots + n_s$ for the unique choice of
the $n_i$ as in Conjectures 4.6.  This is not a universal phenomenon,
however:  $\dim \maxtopdeg(M)$ may be strictly smaller than $\sum_{i \le
s} n_i$; in fact, the equality $\dim
\maxtopdeg(M) = 0$ (i\.e., a unique maximal top-stable degeneration) is
compatible with an infinite poset of top-stable degenerations of $M$; see
\cite{\constructing}.  Essentially all of the ``slack" left in
Conjectures 4.6 can be filled with examples.

\subhead 4.B. Prototypical examples and proof of Theorem 4.5
\endsubhead

We begin with a number of theoretical remarks which are pivotal for all
our computations.  Suppose that $P/C$ degenerates to
$P/C'$, where, as above, $P = \la e$ is a local module; that is, $e$ is a
vertex of the quiver $Q$, and $C, C'\in \grasstd$.  As in 4.A, we let
$\omega_1,
\dots, \omega_m$ be oriented cycles from  $e$  to $e$ which generate $eJe$
modulo $\stab_{eJe} C = \{a \in eJe
\mid Ca \subseteq C\}$.   From Proposition 4.2 we know that there exist rational
functions $p_i(\tau)$ such that
$$\lim_{\tau \rightarrow \infty} Ce(\tau) = C',$$
where $e(\tau) = e + \sum_{1 \le i \le m} p_i(\tau) \omega_i$ is
identified with an element of $e + eJe \cong \unirad$ whenever $\tau \in
\AA^1$ lies outside the union of the pole sets of the
$p_i(\tau)$.  We will write $\tau$ for either a variable
over $K$ or a scalar as convenience dictates.  Adopting the convention
$\deg p = \deg f - \deg g$ if $f,g$ are nonzero polynomials in $K[\tau]$
with $p = f/g$, and $\deg p = - \infty$ when
$p = 0$, we note that $P/C \cong P/C'$ in case none of the
$p_i(\tau)$ has positive degree. 

Proposition 4.2 provides elementwise access to $C'$ by way of the following
observation:  Fix a $K$-basis $(y_j)$ for $JP$.  For any element $x \in
JP$, there then exist unique rational functions $\rho_j(\tau)$ with the
property that 
$x e(\tau) =  \sum_j \rho_j(\tau) y_j$ 
whenever $\tau \in \AA^1$ lies outside the pole
sets of the $p_i(\tau)$.  Note that, for $x \ne 0$, at least one of the
$\rho_j$ is nonzero since, for the pertinent values of $\tau$, right
multiplication by $e(\tau)$ is an automorphism of $P$.  We define the
$e(\tau)$-{\it degree\/} of a nonzero element $x$ in $JP$ as $\max\{\deg
\rho_j(\tau)\}$.
The $e(\tau)$-degree does not depend on our choice of a basis for
$P$, but only depends on $x$.   Moreover, if $x \ne 0$,  we let
$\widetilde{\rho_j} (1/\tau)$ be the rational function in $1/\tau$ for
which 
$$\widetilde{\rho_j} (1/\tau) = (1 / \tau)^{e(\tau)\text{-}\deg(x)}
\rho_j(\tau).$$ 
By construction, all of the $\widetilde{\rho_j}(1/\tau)$ are defined at
$1/\tau = 0$, and at least one of the values $\widetilde{\rho_j}(0)$ is
nonzero.  Thus 
$$x' = \sum_j \widetilde{\rho_j} (0) y_j$$ is a well-defined element
of $JP$, nonzero whenever $x$ is nonzero.  As we will see next, $x'$
is determined by $x$ up to a factor in $K^*$, irrespective of our choice of
reference basis.

\proclaim{Lemma 4.7}  Keep the above notation.  Then
 $\lim_{\tau \rightarrow \infty}
(Kx) e(\tau) = K x'$ for every nonzero element $x \in JP$. 
In particular, $K x'$ is a one-dimensional subspace of $C'$ whenever $x \in
C \setminus\{0\}$. 
\endproclaim 

\demo{Proof} Consider the curve $\tau \mapsto K \bigl(\sum_j
\widetilde{\rho_j}(1/\tau) y_j\bigr)$ in $\Gr(1,JP)$.  It
is defined at the point $\infty$ of $\PP^1$ and coincides with $\tau
\mapsto K x e(\tau)$ on a dense subvariety of $\PP^1$; indeed, 
$Kx e(\tau) = K (1/\tau)^{e(\tau)\text{-}\deg(x)} x e(\tau) =  K\bigl(
\sum_j \widetilde{\rho_j}(1/\tau) y_j \bigr)$ for all $\tau
\in K^*$ for which the $\rho_j(\tau)$ are defined.  \qed
\enddemo              

We will  --  mostly tacitly  --  make use of this fact in the following. 
In each of the examples, completeness of the given list of
top-stable degenerations of $M$ can be verified via 
Proposition 4.2(1).  We will carry out this
argument in Example 4.8 and leave it to the reader in the other
instances.  In the first two of the three families of examples of this
section, the posets of top-stable degenerations are finite.  They show
that, for any value of $\m = \dim \autlap.C$, there may be as many as
$2^{\m}$ top-stable degenerations of
$P/C$ or as few as $\m+1$.  These are the
upper and lower bounds conjectured in 4.6(2).

\definition{Example 4.8}  Let $m\ge 1$ and $\la = KQ /I$, where $Q$ is the
following quiver

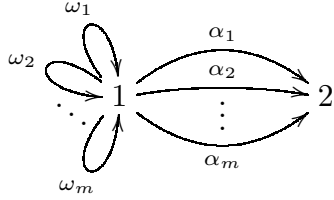
\ignore
$$\xymatrixrowsep{2.0pc}\xymatrixcolsep{2pc}
\xymatrix{ 1 \ar@(ul,u)^{\omega_1} \ar@(ul,l)_{\omega_2}
\save+<-3.5ex,-0.5ex>
\drop{\ddots}\restore \ar@(dl,d)_{\omega_m}
\ar@/^3.5ex/[rr]^{\alpha_1} \ar@/^0.75ex/[rr]^{\alpha_2}
\ar@/_3.5ex/[rr]_{\alpha_m}
 &\save+<0ex,-0.6ex> \drop{\vdots} \restore
&2 }$$
\endignore

\noindent and $I \subseteq KQ$ the ideal generated by
$\omega_i\omega_j$ for all $i,j \in \{1, \dots, m\}$ and
$\alpha_j \omega_i$ for all $i, j$ with $i \ne j$.  For
$T = S_1$,  we consider the  point $C = \sum_{1 \le i \le m}
\la \alpha_i$ in
$\operatorname{\frak{Grass}}^T_{2m+1}$.  Then $\autlap.C \cong \AA^m$,
that is, $m=\m$ in this example.

For each subset $\I \subseteq \{1, \dots, m\}$, the preceding comments
yield the following top-stable degeneration $P/C(\I)$ of $M = P/C$. 
Namely, $C(\I)$ equals
$$\lim_{\tau \rightarrow \infty} C \bigl(e + \sum_{i \in
\I} \tau \omega_i \bigr) = \lim_{\tau \rightarrow \infty} \biggl(
\sum_{i \notin \I} \la
\alpha_i + \sum_{i \in \I} \la \bigl( (1/\tau) \alpha_i + \alpha_i
\omega_i \bigr) \biggr)  = \sum_{i \notin \I} \la \alpha_i + \sum_{i \in
\I} \la \alpha_i \omega_i.$$
\noindent Similarly, one shows that $P/C(\I_1) \degen P/C(\I_2)$ is
equivalent to $\I_1 \subseteq \I_2$.  

To see that all proper top-stable degenerations of $M$ are among the
$P/C(\I)$ with $\I \ne \varnothing$, let $P/C'$ be any such degeneration,
and let
$e(\tau) = e + \sum_i p_i(\tau) \omega_i$ be such that $C' =
\lim_{\tau \rightarrow \infty} C e(\tau)$.  Then the maximum of the degrees
of the $p_i(\tau)$ is positive, and for
$\I = \{i \le m \mid \deg p_i(\tau) >0\}$, one verifies $C' \cong
C(\I)$.      

Therefore, there
are precisely $2^m$ distinct top-stable degenerations of $M$ in
total, including
$M = P/C(\varnothing)$. In fact, for each $h \le m$, the number of
isomorphism types of top-stable degenerations of height $h$ equals $\m
\choose h$ as postulated in 4.6(2); in particular, there is a unique
maximal top-stable degeneration of $M$.   For
$\m = 2$, the full poset can be visualized as follows.  Here a dotted
arrow from
$P/C(\I_1)$ to
$P/C(\I_2)$ indicates that the latter is a minimal degeneration of the
former.  For our graphing conventions, we refer to \cite{\dom} and
\cite{\menace}. 

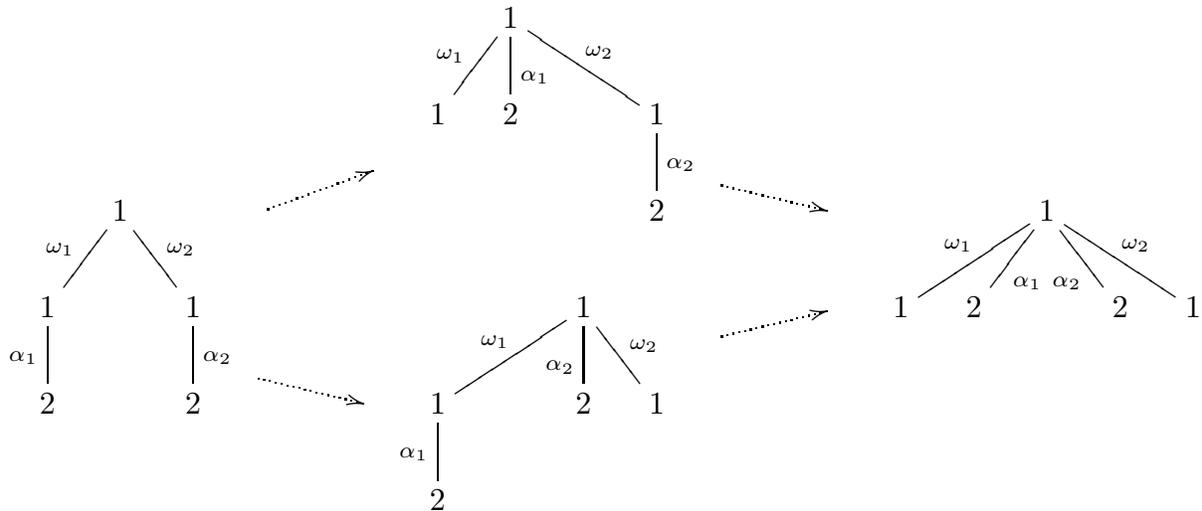
\ignore
$$\xymatrixrowsep{1.5pc}\xymatrixcolsep{1pc}
\xymatrix{
 &&&&& &&1 \edge[dl]_{\omega_1} \edge[d]^(0.6){\alpha_1} \edge[drr]^{\omega_2} \\
 &&&&&\save+<0ex,-4ex> \ar@{<.}[dll] \restore &1 &2 &&1 \edge[d]^{\alpha_2} & \\
 &1 \edge[dl]_{\omega_1} \edge[dr]^{\omega_2} && &&& &&&2 &\save+<0ex,2ex>
\dttdar[rr] \restore && &&&1 \edge[dll]_{\omega_1} \edge[dl]^(0.6){\alpha_1}
\edge[dr]_(0.6){\alpha_2} \edge[drr]^{\omega_2} \\  
1 \edge[d]_{\alpha_1} &&1 \edge[d]^{\alpha_2} & &&&&&1 \edge[dll]_{\omega_1}
\edge[d]_(0.6){\alpha_2} \edge[dr]^{\omega_2} &&\save+<0ex,-2ex> \dttdar[rr]
\restore &&&1 &2 &&2 &1  \\ 
2 &&2 &\save+<0ex,2ex> \dttdar[rr] \restore &&&1 \edge[d]_{\alpha_1} &&2 &1 & &&& \\
 &&&&& &2 }$$
\endignore
\enddefinition

Next, we realize the lower bound on the number of top-stable
degenerations of $P/C$ given in Conjecture 4.6(2) for any value
of $\m = \dim \autlap.C$.  As is predicted by the conjecture, every
number $h$ between $1$ and $\m$ arises as the height of precisely one
top-stable degeneration in this situation. The case $\m=1$ is covered by
the previous example.

\definition{Example 4.9} Let $m \ge 2$ and $\la = KQ/I$, where $Q$ is
the quiver

\ignore
$$\xymatrixrowsep{2.0pc}\xymatrixcolsep{4pc}
\xymatrix{ 1 \ar@(ul,dl)_{\omega} \ar[r]^{\alpha} &2 }$$
\endignore

\noindent and $I = \langle \omega^{m+1} \rangle$.  Once more, we set $T =
S_1$.  We focus on the module $M = P/C$, where $C$ is the point
$\sum_{0 \le i \le m-1} \la \alpha \omega^i$ in
$\operatorname{\frak{Grass}}^T_{m+2}$ and observe that again $\autlap.C
\cong \AA^m$.  Using the above remarks (see \cite{\constructing} for
detail), one checks that the top-stable degenerations of
$M = M^{(0)}$ are precisely the local modules $M^{(j)} = P/C^{(j)}$, where
$C^{(j)} = 
\sum_{0 \le i \le m, i \ne m - j} \la
\alpha
\omega^i$ for
$1 \le j \le m$, as depicted below.

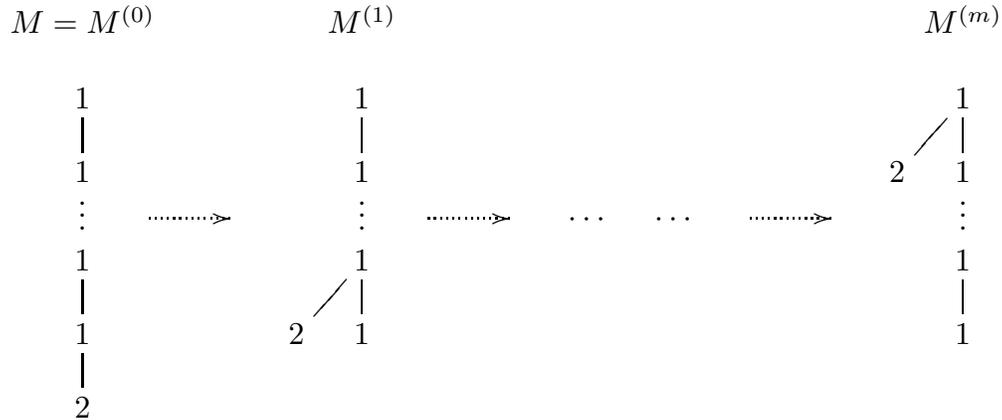
\ignore
$$\xymatrixrowsep{0.2pc}\xymatrixcolsep{0.8pc}
\xymatrix{
 &&& && &&& && &&& &&  \\   \\ 
1\save+<0ex,6ex> \drop{M=M^{(0)}} \restore
\edge[dd] &&& &&1\save+<0ex,6ex> \drop{M^{(1)}} \restore \edge[dd] &&& &&
&&& &&1\save+<0ex,6ex> \drop{M^{(m)}} \restore
\edge[ddl] \edge[dd]  \\   \\  
1 &&& &&1 &&& && &&& &2 &1  \\
\save+<0ex,0.9ex> \drop{\vdots} \restore &\dttdar[rr] &&
&&\save+<0ex,0.9ex> \drop{\vdots} \restore &\dttdar[rr] && &\cdots
&\cdots &\dttdar[rr] && &&\save+<0ex,0.9ex> \drop{\vdots} \restore  \\ 
1 \edge[dd] &&& &&1 \edge[ddl] \edge[dd] &&& && &&& &&1 \edge[dd]  \\   \\ 
1 \edge[dd] &&& &2 &1 &&& && &&& &&1  \\   \\ 
2  }$$

\endignore
\enddefinition

The third class of examples illustrates the infinite scenario and proves
Theorem 4.5.

\definition{Examples 4.10 / Proof of Theorem 4.5} As in Theorem 4.5, let
$l, n_1,\dots,n_s$ be positive integers with $l \ge s$.  In part (a),
we specialize to the case
$s = l = 1$ to construct the basic building blocks for the
general case.  In part (b), we subsequently duplicate and assemble them to
meet the initial requirements of the theorem, combined with variable
saturated chain lengths of top-stable degenerations as postulated in (i)
for $l=s$.  In part (c) we first address the case $s = l$, combined with
constant saturated chain length as postulated in (ii), and in (d)
finally, we indicate how to realize (i) and (ii) for $l > s$.

(a)  We return to the algebra
$\la$ of Example 4.8 with $m = n_1 + 1$ and again let $T = S_1$.  But
this time, we consider the point $C = \sum_{1 \le i \le m} \la
\alpha_i + \la \bigl(\sum_{1 \le i \le m} \alpha_i \omega_i \bigr)$ in
$\operatorname{\frak{Grass}}^T_{2m}$. 

The situation is akin to that of 4.8 in that $\autlap.C \cong \AA^m$, so
again $\m=m$, and each subset $\I \subseteq
\{1, \dots, m\}$ leads to a family of degenerations of $M = P/C$. Once
more, all top-stable degenerations of $M$ are among the ones we display in
the following.  Given $\I$ and a vector $k = (k_1, \dots, k_m) \in K^m$
with
$k_i \ne 0$ precisely when $i \in \I$, we find:
$$\align C(\I,k): &= \lim_{\tau \rightarrow \infty} C \bigl(e + \sum_{1
\le i \le m} k_i \tau \omega_i\bigr) \\ &=  \lim_{\tau \rightarrow \infty}
\biggl(\sum_{i \notin \I} \la
\alpha_i + \sum_{i \in \I} \la \bigl((1/ \tau) \alpha_i + k_i \alpha_i
\omega_i \bigr) + \la \bigl( \sum_{1 \le i \le m} \alpha_i \omega_i
\bigr) \biggr)\\ &=  \lim_{\tau \rightarrow \infty} \biggl(\sum_{i \notin
\I} \la
\alpha_i + \sum_{i \in \I} \la \bigl((1/ \tau) \alpha_i + k_i \alpha_i
\omega_i \bigr) + \la \biggl( \sum_{i \notin \I} \alpha_i \omega_i
 - \sum_{i \in \I} (1/\tau k_i) \alpha_i \biggr) \biggr).
\endalign$$  If $\I$ is a proper subset of $\{1, \dots, m\}$, we  derive 
$$C(\I): = C(\I,k) = \sum_{i \notin \I} \la \alpha_i + \sum_{i \in \I}
\la \alpha_i
\omega_i + \la \bigl( \sum_{i \notin \I} \alpha_i \omega_i \bigr),$$
regardless of the choice of $k_i \in K$.  In particular, the case where
$\{1,
\dots, m\} \setminus \I$ is a singleton $\{j\}$ yields a maximal
top-stable degeneration $P/C(\I)$, since $P/C(\I)$ has Loewy length $2$
 in this situation;  indeed, $C(\I) = \la \alpha _j + \sum_{1
\le i \le m} \la \alpha_i \omega_i$.  

If, on the other hand,
$\I = \{1, \dots, m\}$, the last summand in the final line of the above
limit computation equals
$\la
\bigl(\sum_{i \in \I} (1/k_i) \alpha_i \bigr)$ for $\tau \ne 0$, and
consequently $C(\I,k)$ takes on the form 
$$C(\I,k) = \sum_{1 \le i \le m} \la \alpha_i \omega_i + \la
\bigl(\sum_{1 \le i \le m} (1/k_i) \alpha_i \bigr).$$ Clearly, each such
$P/C(\I,k)$ is again a maximal top-stable degeneration of $P/C$, its
isomorphism type depending only on the point
$\overline{(1/k_1, \dots, 1/k_m)} \in \PP^{m-1}$.  These latter
degenerations thus form a family  parametrized by the torus orbit of
dimension $m - 1$ in $\PP^{m-1}$. 

As for the degenerations $P/C(\I)$, where the cardinality of $\{1,
\dots, m\} \setminus \I$ exceeds $1$:  They are non-maximal, and their
top-stable degenerations can in turn be determined along the above lines,
with the set $\I' = \{1, \dots, m\} \setminus \I$ taking over the role
played by $\I$ in the previous paragraphs.  In particular:  Whenever
$|\{1,
\dots, m\} \setminus \I_1| \ge 2$, the relation ``$P/C(\I_1)$ $\degen$
$P/C(\I_2)$" is equivalent to ``$\I_1 \subseteq \I_2$".  Thus, each set
$\I$ with $|\I'| \ge 2$ leads to another infinite family of maximal
top-stable degenerations $P/C(\I, k)$ with $k_i$ nonzero precisely when
$i \in \I'$.  It is accordingly parametrized by the points of a torus
orbit of dimension $|\I'| - 1$ in
$\PP^{m-1}$; namely,
$C(\I, k) = \sum_{1 \le i \le m} \la \alpha_i \omega_i + \la
\bigl(\sum_{i \in \I'} (1/k_i) \alpha_i \bigr)$.     

 We conclude that the non-maximal top-stable degenerations of
$M$ are precisely the $P/C(\I)$ for which the cardinality of $\{1, \dots,
m\}
\setminus
\I$ is at least
$2$, while the maximal top-stable degenerations of $M$ form a
$\PP^{m - 1}$-family via
$$\overline{(a_1, \dots, a_m)} \mapsto \sum_{1 \le i \le m} \la \alpha_i
\omega_i + \la \bigl( \sum_{1 \le i \le m} a_i \alpha_i \bigr).$$ This
latter family breaks up into subfamilies indexed by torus orbits of
dimension $m - h$, consisting of degenerations of fixed height
$h$, respectively.   Moreover, we glean:  The minimum of the heights of
the maximal top-stable degenerations of $M$ is
$1$, and the maximum is
$m - 1$.

For $m = 3$, the hierachy of top-stable degenerations of $M$ thus takes
on the following form:

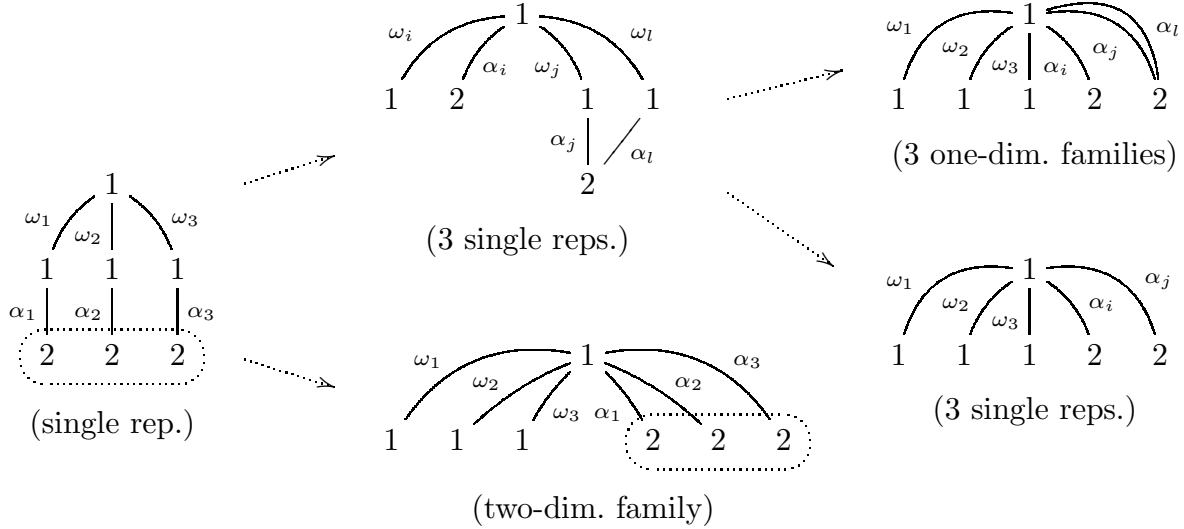
\ignore
$$\xymatrixrowsep{1.2pc}\xymatrixcolsep{0.8pc}
\xymatrix{
 && &&& &&&1 \edge@/_2ex/[dll]_(0.65){\omega_i}
\edge@/_1ex/[dl]^(0.6){\alpha_i} \edge@/^1ex/[dr]_(0.6){\omega_j}
\edge@/^2ex/[drr]^(0.65){\omega_l} && &&&\save+<0ex,-4ex> \ar@{<.}[dll]
\restore  &&&1 \edge@/_3ex/[dll]_(0.65){\omega_1} \edge@/_1ex/[dl]_(0.6){\omega_2}
\edge[d]_(0.6){\omega_3} \edge@/^1ex/[dr]_(0.6){\alpha_i}
\edge@/^3ex/[drr]_(0.6){\alpha_j} \edge@/^4ex/[drr]^(0.7){\alpha_l} \\
 && &&&\save+<0ex,-4ex> \ar@{<.}[dll] \restore  &1 &2 &&1 \edge[d]_{\alpha_j} &1
\edge[dl]^{\alpha_l} & && &1 &1 &1 \dropdown{4}{\txt{ (3 one-dim.~families)}} &2 &2
\\
 &1 \edge@/_1ex/[dl]_(0.6){\omega_1} \edge[d]_(0.6){\omega_2}
\edge@/^1ex/[dr]^(0.6){\omega_3} & &&&
&&&\dropdown{4}{\txt{ (3 single reps.)}} &2 & &\dttdar[drr] && & \\
1 \edge[d]_{\alpha_1} &1 \edge[d]_{\alpha_2} &1 \edge[d]^{\alpha_3} & &&&&& &&&
&&&&&1 \edge@/_3ex/[dll]_(0.65){\omega_1} \edge@/_1ex/[dl]_(0.6){\omega_2}
\edge[d]_(0.6){\omega_3} \edge@/^1ex/[dr]^(0.6){\alpha_i}
\edge@/^3ex/[drr]^(0.65){\alpha_j} \\  
2 \levelpool{2} &2 \dropdown{5}{\txt{(single rep.)}} &2 &&& &&&&1
\edge@/_3ex/[dlll]_(0.65){\omega_1} \edge@/_0.75ex/[dll]_(0.6){\omega_2}
\edge@/_0.5ex/[dl]^(0.6){\omega_3} \edge@/^0.5ex/[dr]_(0.6){\alpha_1}
\edge@/^0.75ex/[drr]^(0.6){\alpha_2} \edge@/^3ex/[drrr]^(0.65){\alpha_3} && &&&1 &1
&1 \dropdown{4}{\txt{ (3 single reps.)}} &2 &2 \\
 && &&&\save+<0ex,4ex> \ar@{<.}[ull] \restore &1 &1 &1 &\dropdown{5}{\txt{
(two-dim.~family)}} &2 \levelpool{2} &2 &2
 }$$
 \endignore
\vskip0.3truein

\noindent Here the dotted enclosure in the graph of $M$ indicates that
the $C$-residues of the paths $\alpha_1 \omega_1$, $\alpha_2 \omega_2$,
$\alpha_3 \omega_3$ are linearly dependent, whereas any two of them are
linearly independent.  The graph in the lower-most row is to be
interpreted analogously;  but, while for $M$ the linear dependence
relationship is prescribed, namely $\sum_{1 \le i \le 3}
\alpha_i
\omega_i = 0$, in the lower-most graph, it takes on the form $\sum_{1
\le i \le 3} a_i \alpha_i = 0$ for arbitrary $(a_1,a_2,a_3) \in (K^*)^3$,
to the effect that this latter graph represents a family indexed by the
dense torus orbit in
$\PP^2$ (a formal description was given above).  Similarly, the top
graph in the last column stands for a family of modules with $a_j
\alpha_j + a_l \alpha_l = 0$, where the pair
$(a_j, a_l)$ traces the homogenoeous coordinates of a $1$-dimensional
torus orbit in $\PP^2$, namely the orbit of the element with homogeneous
coordinates $1$ in positions $j,l$ and $0$ in position $i$.  Throughout,
the indices $i,j,l$ trace all choices yielding $\{i,j,l\} = \{1,2,3\}$.
 
\medskip

(b) Let $n_1, \dots, n_s$ be given, and let $Q$ be the quiver with
$s+1$ vertices, labeled $1$ and $a_1, \dots, a_s$ for $1
\le i \le s$.  Moreover, suppose $Q$ has $\sum_{1 \le i \le s} (n_i+1)$
loops
$\omega_{ij}$ at $1$, where $1 \le i \le s$,
$1 \le j \le n_i +1$, next to another set of $\sum_{1 \le i \le s}
(n_i+1)$ arrows, 
$\alpha_{ij}: 1 \rightarrow a_i$, one for each eligible choice of $i$ and
$j$.  The ideal
$I \subseteq KQ$ is generated by 
$\omega_{ij} \omega_{k\ell}$ for all legitimate choices of $i,j,k,\ell$ and
by the paths $\alpha_{ij} \omega_{k\ell}$ for $(i,j) \ne (k,\ell)$.  Again $P
=
\la e_1$, but this time, 
$$C  = \sum_{\text{all\ } i,j} \la \alpha_{ij} + \sum_{1 \le i
\le s} \la \bigl( \sum_{ 1 \le j \le n_i + 1} \alpha_{ij}
\omega_{ij} \bigr) \ \subseteq \ JP.$$
Mimicking the argumentation under (a), one checks that the maximal
top-stable degenerations of $M = P/C$ form a 
$\PP^{n_1} \times \cdots \times \PP^{n_s}$-family via
$$\biggl( \overline{(k_{11}, \dots, k_{1,n_1 +1})}, \dots,
\overline{(k_{s1},
\dots, k_{s,n_s+1})} \biggr) \mapsto \biggl( \sum_{i,j} \la \alpha_{ij}
\omega_{ij}\  +  \sum_{1 \le i \le s} \la \bigl(  \sum_{1 \le j
\le n_i + 1} k_{ij} \alpha_{ij} \bigr) \biggr).$$ 
Moreover, the minimum
of the heights of the maximal top-stable degenerations of
$M$ is $s$, and the heights of the maximal top-stable
degenerations trace all values between $l = s$ and $\sum_{1 \le i
\le s} n_i$.
\medskip

(c)  To realize examples analogous to those under (b), but satisfying the
additional condition that all maximal top-stable degenerations of $M$
have the same height $l = s$, we once more modify the quiver $Q$.  The
vertices are as under (b), and again we suppose that there are $\sum_{1
\le i \le s} (n_i+1)$ distinct loops $\omega_{ij}$ at $1$, where $1 \le i
\le s$ and $1 \le j \le n_i + 1$; but now we place precisely one arrow
from $1$ to each of the other vertices, say $\alpha_i: 1 \rightarrow
a_i$.  We factor the following ideal $I$ out of the path algebra $KQ$ to
obtain
$\la$: namely $I$ is generated by the products
$\omega_{ij} \omega_{kl}$ for arbitrary choices of $i,j$ and $k, l$, and
by the paths $\alpha_l \omega_{ij}$ for $l \ne i$.  Moreover, we let $M =
P/C$, where $C = \sum_{1 \le i \le s} \la
\alpha_i$.  By means of Proposition 4.7, it is straightforward to check 
that the dimension of $\autlap.C$ is $\m =
\bigl(\sum_{1 \le i \le s} n_i \bigr) + s$.  Along the previous pattern,
one moreover shows that the following family
$\bigl(P/C_{k}
\bigr)$, parametrized by the points $k \in \PP^{n_1}
\times\cdots \times \PP^{n_s}$, gives the maximal top-stable
degenerations of $M$, each isomorphism type arising precisely once: Let
$k = \bigl(\overline{k^{(1)}},
\dots, \overline{k^{(s)}}\bigr)$ be a point in $\PP^{n_1} \times \cdots
\times
\PP^{n_s}$, and let
$\bigl( k_{i 1}, \dots, k_{i, n_{i}+1} \bigr)$ be homogeneous coordinates
for $\overline{k^{(i)}}$; then $C_k \subseteq P = \la e_1$ is generated by
the sums $\sum_{1 \le j \le n_i +1} k_{ij} \alpha_i
\omega_{ij}$ for $1 \le i \le s$.  It is obvious that all of the
degenerations $P/C_k$ of $M$ have height $s$.
\medskip

(d)  For $l > s$, we use the basic setup of (b) or (c)
and, in each case, tack on
$l - s$ replicas of Example 4.8 for $m = 1$.  This means adding on $l-s$
new vertices, $l-s$ new loops at
$1$, and one arrow from $1$ to each of the new vertices, next to
additional relations modeled on 4.8 and relations ensuring that the
new building blocks interfere neither with each other nor with the old
picture.  The modification brings the minimal saturated chain length up to
$l$ and the orbit dimension $\dim \autlap.C$ to
$l + \sum_{1 \le i \le s} n_i$, yet adds no new invariants $n_k$ if $s \ge
1$.  Instead, it effects a parallel-shift of the top-stable degeneration
posets under (b) and (c).
\qed    
\enddefinition

\head 5.  Degenerations of modules with squarefree tops
\endhead

From now on, we let $T \in \lamod$ be a squarefree semisimple module of
dimension $\t$, say
$$T = \la e_1/Je_1 \oplus \cdots \oplus \la e_\t/J e_\t \quad \text{and}
\quad P = \la e_1 \oplus \cdots \oplus \la e_\t$$ 
after renumbering of the vertices of $Q$.  Note that $P$ is a left ideal of
$\la$.  Clearly, the group
$\aut_{\la} (T)$, which arises as a semidirect factor of $\autlap$ (see
Observation 3.2), is a $\t$-dimensional torus in the present situation.  We
fix the following incarnation $\T$ of this torus inside $\autlap$ which, a
priori, is only unique up to conjugation in the full automorphism group: 
namely, we let $\T$ consist of the automorphisms of $P$ sending $z$ to $z
\bigl( \sum_{1 \le i
\le \t} a_i e_i \bigr)$ where $(a_1, \dots, a_{\t}) \in (K^*)^{\t}$; any
such automorphism will be identified with the tuple $(a_1, \dots, a_\t)$. 
Then
$\autlap = \unirad \rtimes \T$.

\subhead 5.A. The basic theorems for squarefree $T$ \endsubhead

In contrast to the local case, there may be a plethora of layer-stable
degenerations of modules with non-simple squarefree top.  In fact, for any
$r \ge 0$, one can find a module $M = P/C$ having a $\PP^r$-family of
such degenerations (see Example 5.9), and there is no uniform upper bound
on the lengths of chains of indecomposable layer-stable degenerations.  

\proclaim{Theorem 5.1. Top-stable degenerations}  Let $T$ be a direct sum
of
$\t$ pairwise nonisomorphic simple modules, $P = \bigoplus_{1 \le i \le
\t} \la e_i$ its projective cover, $d$ a positive integer, and 
$M = P/C$ with $C \in \grasstd$. 
\smallskip 

\noindent {\rm\bf (1) Structure of the $\autlap$-orbits and chain
lengths.} The
$\autlap$-orbits of $\grasstd$  are isomorphic to direct products of full
affine spaces and tori.  More precisely,  
$$\autlap.C  \cong  \AA^{\m} \times (K^*)^{\t - \s},$$  
where $\s = \s(M)$ is the number of indecomposable summands of $M$, and 
$$\m = \m(M) = \dim_K \hom_\la(P,JM) -  \dim_K \hom_\la(M, JM).$$
In particular, the lengths of chains of top-stable degenerations of $M$
are bounded above by $\t - \s$ plus the sum of the multiplicities of the
simple summands of $T$ as composition factors of $JM/ \soc(JM)$.   
\smallskip 

\noindent {\rm\bf (2) Existence.}  The following
conditions are equivalent:
\roster
\item"(i)" $M$ has no proper top-stable degenerations.

\item"(ii)" $\autlap.C$ is a singleton, i\.e., $\m = \t - \s = 0$.

\item"(iii)"  $M$ is a direct sum of local modules, and  every
homomorphism in $\Hom_\la(P,JP)$ leaves $C$ invariant.  {\rm(}The latter
condition is equivalent to requiring that $Cu \subseteq C$  for every
nontrivial path $u$ starting in one of the vertices $e_1, \dots,
e_\t$ and ending in one of these vertices.{\rm)}

\item"(iv)"  $C$ is invariant under all endomorphisms of $P$.
\endroster

\noindent {\rm\bf (3) Unique existence.} $M$ has a unique proper top-stable
degeneration if and only if
$\autlap.C \cong \AA^1$, that is, precisely when $\m = 1$ and $\t - \s =
0$.  In particular, this condition forces $M$ to be a direct sum of local
modules.

If $\m = 0$ and $\t - \s = 1$, then $M$ has precisely two distinct proper
top-stable degenerations.  
\smallskip 

\noindent {\rm\bf (4) Bases.}  If $M'$ is a top-stable degeneration
of $M$, then $M$ and $M'$ share a basis consisting of paths.  More
precisely, there exist submodules $D$ and $D'$ of $P$ with $M \cong P/D$ and
$M' \cong P/D'$, together with a set
$\B$ of paths in $KQ$ {\rm(}which is closed under right subpaths{\rm)}, such
that
$\{p + D \mid p \in \B \}$ is a basis for $P/D$ and $\{p + D' \mid p \in
\B \}$ is a basis for $P/D'$. 
\smallskip

\noindent {\rm\bf (5) The maximal top-stable degenerations of $M$} always
possess a fine moduli space classifying them up to isomorphism, namely 
$$\maxtopdeg(M)  = \{C' \in \overline{\autlap.C} \mid \autlap.C' \text{\
is a singleton}\},$$ 
and the latter is a projective variety of dimension at most $\max\{0, \m +
(\t -
\s) - 1\}$.
\smallskip

\noindent {\rm\bf (6) In case $\m = 0$, the degeneration order equals
$\le_{\kleinext}$.}  Suppose $\m = 0$, and let $M = \bigoplus_{1 \le k
\le \s} M_k$ be a decomposition of $M$ into indecomposable summands.  Then
every top-stable degeneration of $M$ is a direct sum of top-stable
degenerations of the $M_k$.  Moreover, the poset of top-stable
degenerations of $M$ is finite, and $M \degen M'$ if and only if $M
\le_{\kleinext} M'$, for all $d$-dimensional modules $M'$ with top $T$.  

In particular: If $\Hom_\la(P,JP)
= \Hom_\la(P,\soc(P))$, the degeneration order on the class of left
$\la$-modules with top $T$ coincides with the $\Ext$-order.
\endproclaim

On the side, we point out that the situation where $M$ has precisely two
proper top-stable degenerations is not restricted to $\m = 0$ and $\t - \s
= 1$.  It may also occur when $\m = 2$ and $\t - \s = 0$ (see Example 4.9
for $\m = 2$).
Further we note that, in the nonlocal case, the fine moduli
space for the maximal top-stable degenerations, 
$\maxtopdeg(P/C)$, is usually reducible.  Both of our
examples in Section 5.B attest to this.  

We smooth the road towards a proof of the first part of Theorem 5.1. 
A proof of the
first lemma can be found in \cite{\class, Proposition 2.9(3)} 

 \proclaim{Lemma 5.2}   For any $C \in
\grasstd$, there exists a point $C'
\in \autlap.C$ such that the $\T$-orbit $\T.C'$ is isomorphic to the torus
$(K^*)^{\t - \s}$, where $\s = \s(P/C)$.  In fact, this is true whenever
$P =
\bigoplus_{1 \le i \le \s} P_i$ and $C = \bigoplus_{1 \le i \le \s} C_i$
with $C_i \subseteq P_i$ such that $\bigoplus_{1 \le i \le \s}
P_i/C_i$ is a decomposition of $P/C$ into indecomposable summands. \qed
\endproclaim

\proclaim{Lemma 5.3}  If $M \in \lamod$ has squarefree top and $\s=
\s(M)$, then 
$$\dim_K \End_\la(M) = \s \ +\ \dim_K \Hom_\la(M,JM).$$ 
\endproclaim

\demo{Proof}  Suppose $M = \bigoplus_{1 \le i \le \s} M_i$ is a
decomposition of $M$ into indecomposable summands.  Since the top of
$M_i$ does not share any simples with that of $M_j$ for $i \ne j$ by
hypothesis, $\End_\la(M)/\Hom_\la(M,JM)$ $\cong$ $\prod_i
\bigl(\End_\la(M_i)/\Hom_\la(M_i,JM_i) \bigr)$.  We may thus focus on the
situation where $M$ is indecomposable.

Set $T = M/JM$.  Moreover, let $\E$ denote the $\la$-endomorphism ring of
$M$ and $\I$ the ideal $\Hom_\la(M,JM)$.  Then $\I$ is nilpotent and the
factor
$\E/ \I$ embeds canonically into $\End_\la(T)$ as a $K$-algebra.  Since
$T$ is squarefree, the algebra $\End_\la(T)$ is a direct product of
copies of $K$, and so is $\E/\I$; say $\E/\I$ is a direct product of $s$
copies of $K$.  If $s$ were larger than $1$, there would be a nontrivial
idempotent in $\E/\I$, which would lift to a nontrivial idempotent of
$\E$.  But this is ruled out by our indecomposability assumption. 
\qed
\enddemo

\demo{Proof of Theorem 5.1, part (1)} 
 (1)  In light of Lemma 5.2, we may
assume without loss of generality that $\T.C \cong (K^*)^{\t -
\s}$.  Consequently, $\dim \stab_{\T} C$ $=$ $\s$.  

Setting $\unirad = \U$, we start by showing  $\stab_{\autlap} C =
\stab_{\U} C \rtimes \stab_{\T} C$.  For that purpose, we abbreviate
$\autlap$ to $\A$,
$\stab_{\autlap} C$ to $\A_0$, $\stab_{\U} C$ to
$\U_0$, and $\stab_{\T} C$ to $\T_0$.  Clearly, $\U_0 \rtimes \T_0$ is a
subgroup of $\A_0$.  Moreover, the algebraic group $\A_0$ is connected 
--  we will ascertain this in the next paragraph -- and $\U_0 \rtimes
\T_0$ is closed in $\A_0$, since $\U_0$ and $\T_0$ are closed in
$\U$ and $\T$, respectively (see, e\.g\., \cite{\Bor, Chapter I, 1.7}). 
Therefore the desired equality $\A_0 = \U_0 \rtimes \T_0$ is equivalent
to equality of the corresponding dimensions $\dim \A_0$ and $\dim \U_0
\rtimes \T_0$. But this dimension equality is an immediate consequence of
Observation 3.2:  Indeed, in view of $\A.C \cong \A/\A_0$ and $\U.C \cong
\U/\U_0$, we obtain  
$\dim \A_0= \dim_K \End_\la(M) + \dim_K \Hom_\la(P,C)$
and $\dim \U_0 \rtimes \T_0 = \dim_K \Hom_\la(P,C) + \dim_K \Hom_\la(M,
JM) + \s$,
whence Lemma 5.3 fills in what we need.

Thus we only need to back our connectedness claim to secure the desired
stabilizer equality.  For that purpose, it suffices to observe that
the open immersion of $\A$ into the affine space $\hom_\la(P,P)$ restricts
to an open immersion of $\A_0$ into that subspace of $\hom_\la(P,P)$ which
consists of the maps leaving $C$ invariant.  In other words, as
a variety, $\A_0$ is isomorphic to an open subvariety of a linear
$K$-space and consequently connected.  

The group $\T_0$ is connected as well:  Indeed, the obvious open immersion
of $\T$ into $\AA^{\t}$ restricts to an open immersion of $\T_0$ into the
$K$-subspace $\{(a_1, \dots, a_\t) \mid C \bigl( \sum_{1\le i\le \t} a_ie_i
\bigr)
\subseteq C \}$.  Since connected
subgroups of tori are direct factors (see, e\.g\., \cite{\Hum, Ex.6 on
p.108}), we infer that $\T = \T_1 \times \T_0$ for some torus $\T_1$.   

Next, we deduce that the variety
$\A.C \cong \A/\A_0$ is isomorphic to the direct product $\bigl(\U/\U_0
\bigr) \times \bigl( \T / \T_0 \bigr)$.  In light of the splitting of
$\T$, we obtain a projection map  
$$\pi_1: \A = \bigl(\U \cdot \T_1 \bigr) \times \T_0 \rightarrow \U \cdot
\T_1$$  of varieties; it is the geometric quotient of $\A$ by the canonical
right action of $\T_0$ (by \cite{\Bor, Chapter II, Proposition 6.6}, for
instance).  Let $\pi_2: \U \cdot \T_1 \rightarrow (\U \cdot
\T_1)/\U_0$ be the geometric quotient of the subgroup $\U \cdot \T_1$ of
$\A$ by its right
$\U_0$-multiplication, and observe that the map $\pi_2 \circ \pi_1$,
passing to consecutive quotients, has the same fibres as the quotient of
$\A$ by
$\A_0$  --  here again, we let $\A_0$ act by right multiplication on
$\A$.  Consequently,
$\A/\A_0 \cong (\U\cdot \T_1)/\U_0$.  In light of the obvious
variety isomorphisms
$(\U \cdot \T_1)/\U_0 =$ $(\U \rtimes \T_1)/ \U_0 \cong$ $\U/\U_0 \times
\T_1 \cong$ $\U/\U_0 \times \T/ \T_0$, we thus obtain the desired
isomorphism $\A.C \cong \U/\U_0 \times \T/\T_0$.
  
The second factor of this direct product is
$(K^*)^{\t - \s}$.  As for the first:  From Observation 3.2 we obtain
$\U/\U_0 \cong \AA^{\m}$ with $\m$ as postulated.

The final statement of part (1) is clear. \qed
\enddemo

The following consequence of the proof of Theorem 5.1(1) generalizes
Proposition 4.2; it provides another specific instance of the
general approach to finding top-stable degenerations suggested by
Observation 3.7.  We set $e = e_1 + \cdots + e_\t$, whence
$P$ takes on the form $\la e$ and $\autlap$ is identified with the group of
units of $e \la e$.  We start with the prerequisites for describing
the relevant curves in the latter group.  Given $i,j \in \{1, \dots, \t\}$,
not necessarily distinct, we let $\bigl( \omega_{ij}^{(r)} \bigr)_{r
\in R_{ij}}$ denote a collection of paths in $Q$, each from $e_j$ to
$e_i$, such that the
$\omega_{ij}^{(r)}$, with $1 \le i,j \le \t$ and $r \in R_{ij}$, generate
$eJe$ modulo
$\stab_{eJe} C = \{a \in eJe \mid Ca \subseteq C\}$.  In particular,
this means that each $\omega_{ij}^{(r)}$ has positive length and
satisfies $\omega_{ij}^{(r)} = e_i \omega_{ij}^{(r)} e_j$.

\proclaim{Corollary 5.4}  We keep the limit notation
of Proposition {\rm 4.2}.  

{\rm (1)} Suppose that  $P/C$ degenerates to $P/C'$ for some $C' \in
\grasstd$.   Then there exist nonzero polynomials $q_1(\tau), \dots,
q_\t(\tau)$ in $K[\tau]$, next to polynomials $p_{ij}^{(r)}(\tau)$,
such that
$$\lim_{\tau \rightarrow \infty} C e(\tau) = C',$$ 
where 
$$e(\tau) = \sum_{1 \le i \le \t} q_i(\tau) e_i + \sum_{1 \le i,j
\le \t} \ \sum_{r \in R_{ij}} p_{ij}^{(r)}(\tau) \omega_{ij}^{(r)} $$  
is identified with a unit in $e\la e$ whenever $\tau
\in \AA^1$ lies outside the zero sets of the $q_i$.

Moreover, given any flag $\C = (C_{\mu})_{1 \le \mu \le d'}$ of subspaces
$C_\mu$ of $C$,
$$\C' = \bigl( \lim_{\tau \rightarrow \infty} C_1 e(\tau), \dots,
\lim_{\tau \rightarrow \infty} C_{d'} e(\tau) \bigr)$$  
is a flag of subspaces of $C'$.  In
particular, for any subspace $D$ of $C$, the limit $D' = \lim_{\tau
\rightarrow \infty} D e(\tau)$ is a subspace of $C'$
of the same dimension.  If $D$ is a $\la$-submodule of $P$, then so is
$D'$.

{\rm (2)}  Conversely, given polynomials as in {\rm (1)}, let $U
\subseteq \AA^1$ be a dense subset on which all of the $q_i$ are nonzero,
and denote by $\psi$ the extension to
$\PP^1$ of the curve $U \rightarrow \overline{\autlap.C}$ which sends $\tau
\in U$ to $C e(\tau)$.  Then either $\psi$ is constant, or else
$\Im(\psi)$ contains a point in $\overline{\autlap.C} \setminus
\autlap.C$. 
\endproclaim

\demo{Proof}  Without loss of generality, we may assume $C$ to be chosen as
in the proof of Theorem 5.1(1).  Indeed, passage to a different point in
the orbit $\autlap.C$ is harmless, since each such point is of the form
$C\bigl (\sum_{1 \le i \le \t} b_i e_i + \sum_{1 \le i,j \le \t} \sum_{r
\in R_{ij}} a_{ij}^{(r)} \omega_{ij}^{(r)} \bigr)$ for suitable
$a_{ij}^{(r)} \in K$ and $b_i \in K^*$.  We pick a sequence $\epsilon_1,
\dots, \epsilon_{\t-\s}$ in $\T$ such that $\T/\stab_\T(C) = \prod_k K^*
\epsilon_k$.  In view of $\autlap.C \cong \U.C \times \T.C$, we then
obtain an isomorphism of varieties, analogous to that of Lemma 4.1:
$$\biggl(\bigoplus_{i,j,r} K\omega_{ij}^{(r)} \biggr) \times 
\biggl( \prod_k K^* \epsilon_k \biggr) \cong
\autlap.C.$$ 
It maps $\bigl(\omega_{ij}^{(r)}, \epsilon_k)$ to $C\biggl( \bigl(\sum_{1 \le
l \le \t} \epsilon_{kl} e_l \bigr) + \omega_{ij}^{(r)} \biggr)$, if
$\epsilon_k = (\epsilon_{k1},  \dots, \epsilon_{k \t})$ with $\epsilon_{kl}
\in K^*$.  Consequently, any curve
$U \rightarrow \autlap.C$, where $U$ is a dense subset of $\AA^1$, has the
form $C \mapsto Ce(\tau)$ for suitable rational functions
$q_i(\tau)$ and $p_{ij}^{(r)}(\tau)$.  Clearing denominators yields
polynomials $q_i(\tau)$ and $p_{ij}^{(r)}(\tau)$ as required.

For (2), it suffices to note that the variety $\autlap.C$ is quasi-affine
and hence has no complete subvarieties of positive dimension.  \qed
\enddemo

Suppose $P/C$ degenerates to $P/C'$.  In parallel to the approach presented
in Lemma 4.7, we compute individual relations in $C'$ from relations in
$C$ as follows:  Let
$e(\tau)$ be as in Corollary 5.4.  Again, we fix a $K$-basis $(y_k)$ for
$JP$ and define the $e(\tau)$-degree of any element $x \in JP$ as in
Section 4.B.  As before, if $x e(\tau) = \sum_k \rho_k(\tau) y_k$, we let
$\widetilde{\rho_k} (1/\tau)$ be the rational function in $1/\tau$ for
which 
$\widetilde{\rho_k} (1/\tau) = (1 / \tau)^{e(\tau)\text{-}\deg(x)}
\rho_k(\tau)$ and set $x' = \sum_k \widetilde{\rho_k} (0) y_k$.

\proclaim{Lemma 5.5} Let $x$ be a nonzero element of $JP$.  Keeping
the above notation, we obtain:

{\rm (1)}  $\lim_{\tau \rightarrow \infty}
(Kx) e(\tau) = K x'$.  

This is a
one-dimensional subspace of $C'$ whenever $x \in C$. 
\smallskip

\noindent Now suppose that all $p_{ij}^{(r)}(\tau)$ are zero, i.e.,
$e(\tau) = \sum_{1 \le i \le \t} q_i(\tau) e_i$.  Then:

{\rm (2)}  The $e(\tau)$-degree of $x$ equals the maximum of the degrees
$\deg q_i(\tau)$ for which $x e_i \ne 0$.  Moreover,
$$Kx' = K\bigl(\sum_{i \in \L} \widetilde {q_i}(0) x e_i \bigr),$$ 
where $\L = \{i \le \t \mid \deg q_i(\tau) = e(\tau)\text{-}\deg(x)\}$ and
$\widetilde{q_i}(1/\tau)$ is the rational function in
$1/\tau$ for which    
$\widetilde{q_i} (1/\tau) = (1 / \tau)^{e(\tau)\text{-}\deg(x)} q_i(\tau)$.
In particular, the $e(\tau)$-degrees of $x$ and $x'$ coincide.

{\rm (3)}  For every subspace $D \subseteq JP$, there exists a
basis $(d_k)$ such that 
$$\lim_{\tau \rightarrow \infty} D e(\tau) = \bigoplus_k \lim_{\tau
\rightarrow \infty} (Kd_k) e(\tau) = \bigoplus_k Kd'_k.$$
More strongly, if $r$ belongs to $\ZZ\cup\{-\infty\}$ and
$D_r \subseteq D$ is the subspace consisting of all elements in $D$ of
$e(\tau)$-degree at most
$r$, any basis for $D_r$ of the above description can be supplemented to a
basis for $D$ with the same  property. 
\endproclaim

\demo{Proof}  Part (1) is proved as Lemma 4.7.  For (2), we observe that
the nonzero candidates among the elements $x e_i$ are linearly
independent, whence it is harmless to incorporate them in the
reference basis
$(y_k)$ of $JP$; indeed, in light of part (1), a basis change leaves
$x'$ invariant up to a factor in $K^*$.  This makes the first assertion
under (2) obvious.  For
$\tau$ in a dense subset of $K$, the product
$(1/\tau)^{e(\tau)\text{-}\deg(x)} x e(\tau)$ consequently takes on the
form
$\sum_{i,\, xe_i \ne 0} \widetilde{q_i} (1/\tau) x e_i$, whence passage to
the limit ``$\tau \rightarrow \infty$" proves the second assertion under
(2).

(3)  By hypothesis, $e(\tau) = \sum_{1 \le i \le \t} q_i(\tau) e_i$
with polynomials $q_i \in K[\tau]$.   Let
$r_1 < \cdots < r_s$ denote the distinct $e(\tau)$-degrees attained on the
nonzero elements of $D$, and consider the partial flag $D_{-\infty} = 0
\subsetneq D_{r_1} \subsetneq \cdots \subsetneq D_{r_s} = D$, where
each $D_j \subseteq D$ is the subspace consisting of the elements of $D$
that have
$e(\tau)$-degree at most $j$.  For every $j \in \{r_1, \dots,
r_s\}$, we moreover define  $\L_j = \{1 \le i \le \t \mid \deg
q_i(\tau) = j\}$ and $f_j = \sum_{i \in \L_j} e_i$.  Part (2) of the
lemma then implies that 
$D_j = D \bigl(\sum_{i \le j} f_i \bigr)$.
Refine the partial flag $(D_j)_j$ of $D$  to a full
subspace flag $0 = B_0 
\subsetneq B_1 \subsetneq \cdots \subsetneq B_t = D$, and suppose
that, for some $m < t$, the family $(b_k)$ is a basis for $B_m$ with the
following properties:  it contains bases for those subspaces $D_j$ which
are contained in $B_m$, and 
$$B'_m = \lim_{\tau \rightarrow \infty} B_m
e(\tau) = \bigoplus_k \lim_{\tau \rightarrow \infty} (K b_k) e(\tau).$$  
Pick $x \in B_{m+1} \setminus B_m$.   Assuming that $\lim_{\tau \rightarrow
\infty} (Kx) e(\tau)$ is contained in $B'_m$, our choice of the $b_k$
provides us with scalars $a_k
\in K$ such that $x' = \sum_k a_k b'_k$; here $x'$ and $b'_k$ result from
$x$ and $b_k$ as described above.  Set $y = \sum_k a_k b_k \in B_m$,
and let $l = e(\tau)$-$\deg(x)$, $l_k = e(\tau)$-$\deg(b_k)$.  The
construction of our flag guarantees that $l_k \le l$ and that all elements
in $D$ of $e(\tau)$-degree $<l$ belong to $B_m$.  In view of part (2) of
the lemma, $x' = x' f_l$ and $b'_k = b'_k f_{l_k}$, whence linear
independence of the $b'_k$ yields
$a_k = 0$ whenever $l_k < l$.  Since, by construction, $y \ne 0$ and $D_{l-1}
\subset B_m$ has basis $\{b_k \mid e(\tau)-\deg(b_k) < l\}$, we glean that
the $e(\tau)$-degree of $y$ also equals $l$ and that $y' =
\sum_k a_k b'_k = x'$.  In other words, $\sum_{i \in \L_l}
\widetilde{q_i} (0) y e_i$ $=$ $\sum_{i \in \L_l} \widetilde{q_i} (0) x
e_i$.  Since $\widetilde{q_i} (0) \ne 0$ for all $i \in \L_l$, we
infer $y e_i = x e_i$ for $i \in \L_l$.  In light of $x e_i =
y e_i = 0$ whenever $i \in \bigcup_{j > l} \L_j$, we conclude that the
$e(\tau)$-degree of $x - y$ is strictly smaller than $l$; again we invoke
part (2) of the lemma.  But, as we noted above, this places the difference
$x - y$ into $B_m$, which is incompatible with our
choice of $x$ outside $B_m$.  An obvious induction now completes the proof
of (3).     
\qed
\enddemo   

From part (3) of Lemma 5.5 we glean in particular that the ``torus
degenerations" of $P/C$ can be obtained along simplistic curves in
$\autlap.C$, because only certain coefficients of the $q_i(\tau)$ matter,
next to their degree differential.  Indeed, $\lim_{\tau
\rightarrow \infty} C
\bigl(\sum_{1 \le i \le \t} q_i(\tau) e_i \bigr)$ remains unchanged if we
replace the $q_i(\tau)$ by monomials $a_i (1/\tau)^{m_i}$ with suitable
coefficients
$a_i \in K^*$ such that $m_i < m_j$ if and only if $\deg q_i > \deg q_j$.

\proclaim{Lemma 5.6}  Let $\T$ be the torus in $\autlap$, as introduced at
the beginning of Section 5, and let $C \in \grasstd$ be a point with $\dim
\T.C > 0$.  Then $\overline{\T.C} \setminus \T.C$ contains at least two
points belonging to different $\autlap$-orbits.  In particular, $P/C$
has at least two nonisomorphic proper top-stable degenerations.
\endproclaim

\demo{Proof}  By part (1) of Theorem 5.1  --  already established  --  we can
write $M = P/C$ in the form $U \oplus V$, where $U$ is indecomposable and
nonlocal.   Without loss of
generality, 
$C = A
\oplus B$ with $A \subseteq \bigoplus_{1 \le r \le u} \la e_r$ and $B
\subseteq \bigoplus_{u+1 \le r \le \t} \la e_r$ such that $U \cong
 (\bigoplus_{1 \le r \le u} \la e_r)/ A$ and $V = 
(\bigoplus_{u+1 \le r \le \t} \la e_r)/B$ for a suitable index $u \in \{2,
\dots, \t\}$. 

First we consider, for any $\tau \in K^*$, the automorphism $f_\tau \in
\T$ defined by $f_\tau (e_1) =
\tau e_1$ and $f(e_r) = e_r$ for $r \ge 2$.  The curve $K^* \rightarrow
\T.C$, given by $\tau \mapsto f_\tau (C)$, has a unique extension
$\PP^1 \rightarrow \overline{\T.C}$, due to completeness of
$\overline{\T.C}$.  We denote the value of this extension at
infinity by $C' =
\lim_{\tau \rightarrow \infty} f_\tau (C)$.  First we note that $C'$ is
again a $d'$-dimensional subspace of $JP$.  To see that $C'$ does not
belong to $\T.C$, we let
$\pi_1: P \rightarrow \la e_1$ be the projection along
$\bigoplus_{r \ge 2} \la e_r$.  Setting $\mu = \dim \pi_1(C)$, we
pick elements
$a_1, \dots, a_\mu \in A$ such that $\pi_1(a_1), \dots, \pi_1(a_\mu)$
form a basis for $\pi_1(A) = \pi_1(C)$.  We supplement these elements with a
basis $a_{\mu +1},
\dots, a_{\nu}$ for $A \cap \Ker(\pi_1)$ to obtain a basis $a_1, \dots,
a_\nu$ for $A$.  Finally, we add on a basis $b_{\nu +1}, \dots, b_{d'}$
for $B$, which results in a basis 
$a_1, \dots, a_{\nu}, b_{\nu +1}, \dots b_{d'}$ for $C$.
Clearly, $a_{\mu +1}, \dots, a_{\nu}, b_{\nu + 1}, \dots, b_{d'}$ are
fixed by all
$f_\tau$, and hence belong to $C'$.  Moreover, the following
spaces are contained in $C'$ (cf\. Observation 3.7):  Namely,
the one-dimensional subspaces $\lim_{\tau \rightarrow \infty} f_\tau(K
a_r)$ for $r$ between $1$ and $\mu$.  If $a_r = \sum_{1 \le i \le \t}
\lambda_{ri} e_i$ with $\lambda_{ri} \in \la$, then the latter space
equals 
$$\lim_{\tau \rightarrow \infty} K \bigl(\lambda_{r1} e_1
+  \sum_{i \ge 2} (1/\tau) \lambda_{ri} e_i \bigr) = K \lambda_{r1}
e_1 =  K \pi_1(a_r).$$ 
Since the elements $\pi_1(a_1), \dots, \pi_1(a_\mu), a_{\mu +1}, \dots,
a_{\nu}, b_{\nu + 1}, \dots, b_{d'}$ of $C'$ are linearly independent by
construction, they form a basis for $C'$.  Consequently, $C' = A_1 \oplus A_2
\oplus B$, where $A_1 = \pi_1(C) = \pi_1(A)$ and $A_2 = A
\cap (\bigoplus_{2
\le r \le u} \la e_r) = \sum_{\mu + 1 \le r \le \nu} K a_r$.  We thus obtain
the following decomposition of the
$\la$-module
$P/C'$:
$$P/C' = \bigl(\la e_1/ A_1 \bigr) \oplus \biggl( \bigl(\bigoplus_{2 \le
r \le u} \la e_r \bigr) / A_2 \biggr) \oplus \biggl(\bigl(\bigoplus_{u+1
\le r\le d'} \la e_r \bigr)/ B \biggr).$$
This module is not isomorphic to $P/C$, since it has more
indecomposable summands than $P/C$.

Analogously, one obtains a point $C'' = A_1' \oplus A_2' \oplus B \in
\overline{\T.C} \setminus \T.C$ with $A_1' = A \cap \la e_1$ and $A_2' =
\pi'(A)$, where $\pi'$ is the projection of $P$ onto $\bigoplus_{2 \le r \le
u} \la e_r$ along $\la e_1 \oplus \bigoplus_{u+1 \le r \le \t} \la e_r$.  

Since $U$ is indecomposable, $\dim A_1' < \dim A_1$, which shows $P/C'$
and $P/C''$ to be nonisomorphic. \qed \enddemo

\demo{Proof of Theorem 5.1, parts (2)--(6)} 

(2). That (i) implies
(ii) follows from part (1), in view of the fact that closedness of the
quasi-affine variety $\autlap.C$ in the complete variety $\grasstd$
forces the dimension of the former to be zero.  The converse is trivial. 
To see that (ii) is equivalent to (iii), note that the equality $\t
= \s$ just means that $M$ is a direct sum of local modules.  Moreover, $\m
= 0$ is tantamount to the equality $\U.C = \{C\}$ (cf\. Observation 3.2);
the latter in turn says $(\id + f)(C) = C$ for all $f
\in \Hom_\la(P,JP)$, which is the same as $f(C) \subseteq C$ for all $f
\in \Hom_\la(P,JP)$.  The implication (iv)$\implies$(ii) is clear, since
(iv) means $\autlap.C$ is a singleton; for (ii)$\implies$(iv), let
$f \in \End_\la(P)$, and consider the automorphism
$g = a \id - f$ of $P$, where $a \in K$ fails to be an eigenvalue of $f$. 
Invariance of $C$ under $g$ then amounts to invariance of $C$ under
$f$.    
 
(3). By the proof of part (1), we may, without loss of
generality, assume that the variety $\autlap.C$ is isomorphic to the
product $\U.C \times \T.C$, with the second factor isomorphic to
$(K^*)^{\t - \s}$.

Whenever $\dim \T.C > 0$, Lemma 5.6 guarantees the
existence of two distinct proper top-stable degenerations of $M$. 
So, on the assumption that $M$ has a unique proper top-stable degeneration,
$\autlap.C$ equals $\U.C \cong \AA^\m$, and as in the proof of Theorem
4.4, part (3b), we derive $\m = 1$.
That, conversely, $\autlap.C \cong \AA^1$ forces $M$ to have a unique
top-stable degeneration, can also be shown as in the local case.  

To justify the final assertion under (3), suppose $\m = 0$ and $\t - \s =
1$.  The existence of at least two nonisomorphic proper top-stable
degenerations of $M$ follows from Lemma 5.6.  On the other hand, 
$|\overline{\T.C} \setminus \T.C| \le 2$, since the embedding $\T.C \cong
K^* \rightarrow \overline{\T.C}$ extends to a morphism $\PP^1 \rightarrow
\overline{\T.C}$ and the image of the latter is closed (see, e.g.,
\cite{\Bor, Chapter AG, Section 7.4}).  This completes the proof of (3). 

Proofs for parts (4) and (5) carry over verbatim from the local
situation.   

(6). Now suppose that $\m = 0$, meaning that the unipotent radical of
$\autlap$ acts trivially on $C$.  Thus $\autlap.C = \T.C$.  Once again, we 
assume that
$C$ is chosen as in the proof of Theorem 5.1(1).  In other words,   
there exists a partition
$\I_1 \cup \cdots \cup \I_{\s}$ of $\I = \{1, \dots, \t\}$ 
with the following property:  For 
$P_k = \bigoplus_{i \in \I_k} \la e_i$, we obtain a decomposition
$C = \bigoplus_{1 \le k \le \s} C_k$ with $C_k \subseteq P_k$ which, in
turn, gives rise to a decomposition  $\bigoplus_{1 \le k \le \s} P_k/C_k$
of $M$ into indecomposable direct summands. 

We first prove that every top-stable degeneration of $M$ is a direct sum
of top-stable degenerations of the $P_k/C_k$.
Let $C' \in \grasstd$ with $P/C \degen P/C'$.  According to Corollary 5.4, 
we choose
$e(\tau)$ such that $C' = \lim_{\tau \rightarrow \infty} C e(\tau)$.  In
our present scenario, we may take $e(\tau)$ to be of the form
$\sum_{1 \le i \le \t} q_i(\tau) e_i$ for suitable polynomials
$q_i(\tau) \in K[\tau]$, and, in view of the remark following Lemma 5.5, we
may simplify $e(\tau)$ further as follows:  Namely, if the degrees of the
$q_i(\tau)$ take
$r+1$ different values, $m_0 > \cdots > m_r$, we may replace $q_i(\tau)$ by
$a_i (1/\tau)^{\nu}$ if $\deg q_i = m_\nu$, where the $a_i$ are nonzero
scalars determined by the $q_i$.  

Clearly, 
$P_k e(\tau) \subseteq P_k$, whenever $\tau \in K^*$.  Consequently,
$\tau \mapsto C_k e(\tau)$ defines a curve from $K^*$ to
$\Gr(\dim C_k, JP_k)$.  Due to closedness of 
$\Gr(\dim C_k, JP_k)$ in $\Gr(\dim C_k, JP)$, we infer that
$\lim_{\tau \rightarrow \infty} C_k e(\tau)$ 
is a $(\dim C_k)$-dimensional subspace of $J P_k$.  A comparison of
dimensions shows that $C'$ is the direct sum of these limits.  Call them
$C'_k$, respectively.  Therefore, $M'= P/C'$ decomposes in the form
$M' = \bigoplus_{1 \le k \le
\s} P_k/C'_k$ with $P_k/C_k \degen P_k/C'_k$ as required.  

Next we show that $M \le_{\kleinext} M'$.  (Our argument will actually
prove that the following holds, irrespective of the value of
$\m$:  In case $C' \in\overline{\T.C}$, the degeneration $P/C'$ satisfies
$P/C \le_{\kleinext} P/C'$.)  The preceding
paragraph legitimizes restriction to the case where $M$ is
indecomposable.  Moreover, it is harmless to assume
$\t > 1$, for otherwise $\T.C$ is a singleton by part (1).   

The crucial step consists of showing that either $P/C \cong P/C'$, or else
there is a nonempty proper subset $\L$ of $\I = \{1, \dots, \t\}$ such that
the top-stably embedded submodule $U = \sum_{i \in
\L} \la (e_i + C)$ of $P/C$ is isomorphic to a direct summand $U'$ of
$P/C'$, with $M/U$ degenerating to a direct sum complement of $U'$ in
$P/C'$.

The $e(\tau)$-degrees of the nonzero elements of $C$ vary
between $-r$ and $0$ in our present setting.  For each integer $k$ in
this range, we consider the
submodule $D_k$ of $C$ consisting of all elements having
$e(\tau)$-degree at most $k$.  Evidently, $D_0 = C$.  According to Lemma
5.5(3), we can find a basis $B$ of $C$ which includes a basis $B_k$ for
each $D_k$ and has the property that 
$$D_k' := \lim_{\tau \rightarrow \infty} D_k e(\tau) = \bigoplus_{b \in
B_k} \lim_{\tau \rightarrow \infty} (Kb)e(\tau)
= \bigoplus_{b \in B_k} K b'$$ 
for all $k$, where the $b'$ are as in Lemma 5.5.  Note that $D'_k$ is a 
submodule of $C'$ by Corollary 5.4.  We introduce
another partition of the set $\I$:  Namely $\I$
is the disjoint union of the subsets $\L_k = \{1 \le i \le \t \mid \deg
q_i(\tau) = k\}$ for $-r \le k \le 0$.  

Suppose $b \in B_k \setminus B_{k-1}$ for some $k > -r$.  This means $b e_i
\ne 0$ for at least one index $i \in \L_k$ and $be_i = 0$ for all $i \in
\bigcup_{j > k} \L_j$.  To adjust the notation to that of Lemma 5.5(2) for
an inspection of $b'$, we set 
$\widetilde{q_i}(1/\tau) = q_i(\tau)(1/\tau)^{k}$.  Whenever $i \in
\L_k$, this gives $\widetilde{q_i}(0) = a_i$.  Lemma 5.5 thus yields 
$\lim_{\tau \rightarrow \infty}(Kb)e(\tau) = K \bigl(\sum_{i \in \L_k}
a_i b e_i\bigr)$, and hence this latter space is contained in $C' \cap
\bigoplus_{i \in \L_{k}} \la e_i$.  For $k = -r$, our choice of $B$ thus
implies
$D'_{-r}  = \lim_{\tau \rightarrow \infty} D_{-r} e(\tau) \subseteq
C'\cap Q_1$, where $Q_1 =
\bigoplus_{i \in \L_{-r}} \la e_i$.  Moreover, we see that $D'_{-r}$ is
isomorphic to $D_{-r}$ via a torus automorphism of $Q_1$, which guarantees
that $Q_1/D_{-r}$ is isomorphic to $Q_1/D'_{-r}$.  
If we denote by $D$ the subspace of
$C$ spanned by the $b \in B \setminus    B_{-r}$, we similarly see that
$D': = \lim_{\tau \rightarrow \infty} De(\tau)$ is contained in $C' \cap
Q_2$, where $Q_2 = \bigoplus_{i \in \bigcup_{-r < k \le 0} \L_k} \la
e_i$.   Let $f$ be the sum of the idempotents generating
$Q_2$; i.e., right multiplication by $f$ projects $P$ onto $Q_2$ along
$Q_1$.  Then
$Df = Cf$, and we find that this submodule of $JP$ is shifted to $D'$ as
we move it to $\infty$ along our given curve:  Indeed, by construction,
$$D' = \bigoplus_{b \in B \setminus B_{-r}} \lim_{\tau
\rightarrow \infty} (Kb)e(\tau) = \bigoplus_{b \in B\setminus B_{-r}}
\lim_{\tau \rightarrow \infty} (Kb)f e(\tau),$$ 
where the last limit equals $\lim_{\tau \rightarrow \infty} (Cf)e(\tau)$,
due to the dimension equality $\dim D = \dim Cf$. In particular, $D'$ is a
submodule of $C'$, since $Cf$ is a submodule of $C$. This confirms that
$Q_2/Cf$ degenerates to $Q_2/D'$.   

In case $r = 0$, i.e., in case all
$q_i(\tau)$ have the same degree, $C'$ equals $C\bigl( \sum_{1
\le i \le \t} a_i e_i \bigr)$ and consequently belongs to $\autlap.C$, that
is, $P/C
\cong P/C'$.  So suppose that $r > 0$, and set $\L = \L_{-r}$.  
Then the submodule $U = \sum_{i \in
\L} \la (e_i + C)$ of $P/C$ satisfies the requirements spelled out above: 
Indeed, $C' = D'_{-r} \oplus D'$ by our choice of $B$, and therefore
$M' = \bigl(Q_1/D'_{-r}\bigr) \oplus \bigl(Q_2/D'\bigr)$. The first
summand is isomorphic to $Q_1/D_{-r} = Q_1/(Q_1 \cap C)
\cong U$.  For the second, we obtain $M/U \cong \bigl( Q_2/Cf \bigr) \degen
\bigl(Q_2/D'\bigr)$ as explained above.

To deduce that $P/C \le_{\kleinext} P/C'$, it now suffices to
keep in mind that 
$D'$ equals the limit $\lim_{\tau \rightarrow \infty} (Cf) \bigl(\sum_{i
\in \I \setminus \L_{-r}}q_i(\tau) e_i\bigr)$.  This puts $D'$ into
the closure of the torus orbit of $Cf$ in $\Gr(\dim Cf, JQ_2)$, whence an
obvious induction completes the argument.

Finally, we observe that there are only $2^{\t}$ possibilities for the
submodule $U \subseteq M$ described above, whence we conclude that $M$
has only finitely many top-stable degenerations, up to isomorphism. 
\qed\enddemo

\proclaim{Theorem 5.7. Layer-stable degenerations} Let $T$ be a direct sum
of $\t$ pairwise nonisomorphic simple modules, $P = \bigoplus_{1 \le i \le
\t} \la e_i$ its projective cover, $d$ a positive integer, and 
$M = P/C$ with $C \in \grasstd$. Let $\s= \s(M)$ and $\m= \m(M)$.
\smallskip 

\noindent{\rm \bf (1) Existence.}
\smallskip  

{\bf (a)} If $M$ is a direct sum of
local modules, that is if $\t - \s = 0$, then  $M$ has no proper
layer-stable degenerations.
\smallskip

{\bf (b)} If $\m = 0$, the following
statements are equivalent:
\roster
\item"(i)"  $M$ has a proper layer-stable degeneration.

\item"(ii)" One of the indecomposable summands of $M$ has a nonzero proper
layer-stably embedded submodule.
\endroster

\noindent In this case, $M$ has only finitely many proper layer-stable
degenerations, the minimal ones being of the form $U \oplus M/U$, where
$U$ is layer-stably embedded in $M$. 
\smallskip 

{\bf (c)}  For any positive integer
$r$, there exists an indecomposable finite dimensional module $M$ with
squarefree top, over a suitable finite dimensional algebra $\la$, such
that $M$ has a
$\PP^r$-family of pairwise non-isomorphic indecomposable layer-stable
degenerations.
\smallskip

\noindent{\rm \bf (2) Bases.}  If $M' = P/C'$ is a layer-stable
degeneration of
$M$, then $M$ and $M'$ share a skeleton {\rm(cf. Section 2)}.  
\endproclaim

\demo{Proof} (1)(a) is due to the fact that $\T$ acts trivially on $C$,
whenever $M = P/C$ is a direct sum of local modules: This means $\autlap.C
= \U.C$, where $\U$ is the unipotent radical of $\autlap$. 
But, letting $\SS = \SS(M)$ be the radical layering of $M$, we know from
Observation 3.2 that $\U.C$ is closed in $\grassSS$, whence Observation
3.1 allows us to conclude that $M$ has no proper layer-stable degeneration.

In light of Observation 3.3, the assertion (1)(b) is a special case of
part (6) of Theorem 5.1, and (1)(c) will be justified in Example 5.9.

(2).  Let $M' = P/C'$ be a layer-stable degeneration of $P/C$ and $\S$
any skeleton of $M'$.  Layer-stability means that $C'$ belongs to the
closure of $\autlap.C$ in $\grassSS$, where $\SS = \SS(M) = \SS(M')$
(Observation 3.1).  Since $\grassS$ is open in $\grassSS$ (Theorem 2.1),
we infer that $\grassS$ intersects $\autlap.C$ nontrivially.  This means
that $\autlap.C \subseteq \grassS$, because $\grassS$ is stable under the
$\autlap$-action, and thus proves our claim.
\qed 
\enddemo

\subhead{5.B Examples and proof of Theorem 5.7, part (1)(c)} \endsubhead

Our approach to computing the layer-stable degenerations of a module $M
= P/C$ (without computing  all of its top-stable degenerations) is
based on the fact that the closure of
$\autlap.C$ in $\grass{\SS(M)}$ is the union of the closures of the orbit 
$\autlap.C$ in the affine varieties $\grassS$, where $\S$ runs through the
skeletons of
$M$.  In writing $e = e_1 + \cdots + e_\t$, we obtain $\autlap.C$ as
the set of points $Cu$ in $\grasstd$, where $u$ runs through the units in
$e \la e$.  One can pare
down the computation required by noticing that, in fact, $\autlap.C$
consists of the points $Cu$, where $u$ is a unit in $e \la e$ which is a
$K$-linear combination of ($I$-residue classes of) paths in $\S$; this is
due to the fact that every module isomorphic to $M$ again has skeleton
$\S$ and thus results from $M$ via an isomorphism having a lift to
$\autlap$ of the described ilk.  

Example 5.8 illustrates the finite scenario.  In contrast to the situation
$\m = 0$ described in Theorem 5.7, even an
indecomposable module $M$ with only finitely many top-stable degenerations
may have indecomposable candidates among its proper layer-stable
degenerations.  Our first example exhibits this phenomenon.  Diagram 5.7
shows the full poset of proper top-stable degenerations of the considered
module.

\definition{Example 5.8}  Let $Q$ be the following quiver:

\ignore
$$\xymatrixrowsep{2pc}\xymatrixcolsep{3pc}
\xymatrix{ 2 \ar[r]^{\beta} &1 \ar@(ul,ur)^{\alpha} &3 \ar[l]_{\gamma} }$$
\endignore

\noindent Define $\la= KQ/ \langle \alpha^3 \rangle$.  Moreover, let
$T = S_1 \oplus S_2 \oplus S_3$, and $d = 5$.  The module whose top-stable
degenerations we want to find is $M = \la /C$, where $C = \la (\alpha
-\beta) + \la (\alpha -\gamma)$.   Clearly, $M$ has precisely three
distinct skeletons, namely:
$$\xalignat3 \S_1 &= \{ e_1, e_2, e_3, \alpha, \alpha^2 \}, &\S_2 &=
\{ e_1, e_2, e_3, \beta, \alpha\beta \}, &\S_3 &= \{ e_1, e_2, e_3,
\gamma, \alpha\gamma \}. \endxalignat$$
Let us abbreviate $\S_1$ to $\S$ and start by determining polynomials for
$\grassS$.  There are precisely two relevant $\S$-critical pairs, namely
$(\beta,e_2)$ and $(\gamma,e_3)$ (we may ignore the pair $(\alpha,
\alpha^2)$ since $\alpha^3 M = 0$), with
$\S(\beta, e_2) = \{\alpha,
\alpha^2 \} = \S(\gamma,e_3)$, leading to the basic congruences
$$\beta \seq X_1 \alpha + X_2 \alpha^2 \quad \text{and} \quad  \gamma
\seq X_3 \alpha + X_4 \alpha^2.$$   
The procedure employed in Example 2.2 yields $\grassS =
\AA^4$.  Note that $M$ corresponds to the point $C = (1,0,1,0)$ in the
present coordinatization.
Using the comments we made at the beginning of this section, we further
obtain 
$$\autlap.C = \{C' = (c'_1, c'_2, c'_3, c'_4) \mid c'_1, c'_3 \ne 0, \
c'_1 c'_4 = c'_2 c'_3 \}$$    
in the coordinates of $\grassS$.  Thus, $\autlap.C \cong \AA^1 \times
(K^*)^2$ in the format of Theorem 5.1(1).

In computing the closure $\overline{\autlap.C}^{\S}$ of $\autlap.C$ in
$\grassS$, one arrives at the following degenerations of $M$; a graph of
$M$ is displayed on the left:

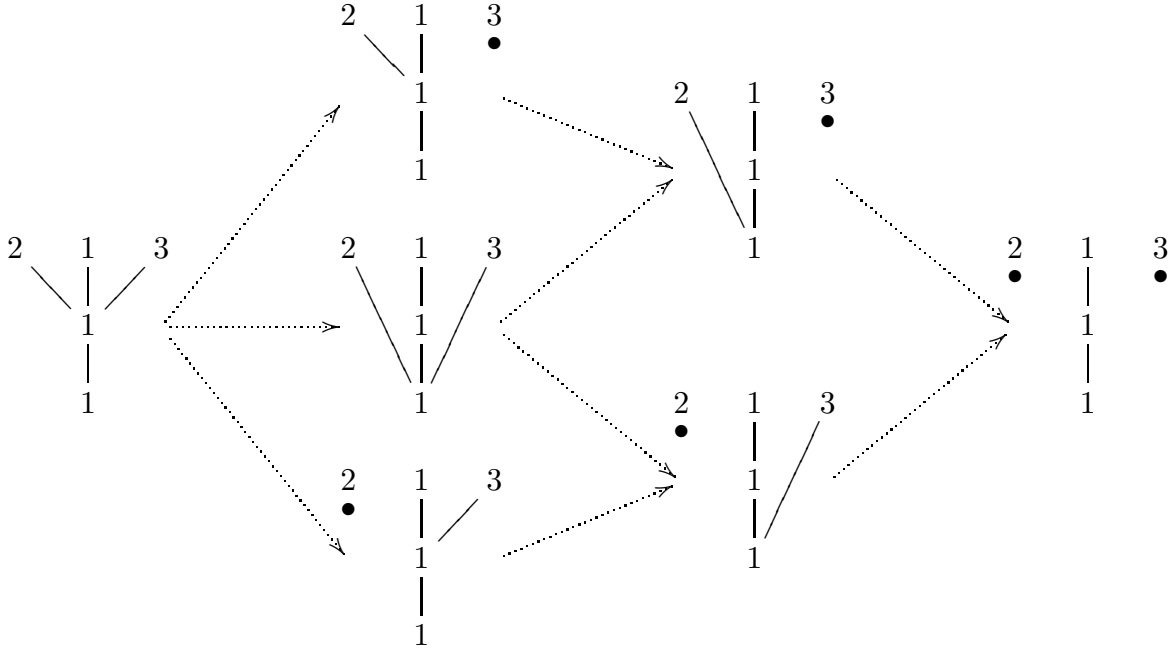
\ignore
$$\xymatrixrowsep{1.0pc}\xymatrixcolsep{1.0pc}
\xymatrix{
 &&&&&2\edge[dr] &1\edge[d] &3\drbl\\
 &&&&&&1\edge[d] &\dttdar[drrr] &&&2\edge[ddr] &1\edge[d] &3\drbl\\
 &&&&&&1 &&&&&1\edge[d] &\dttdar[ddrrr]\\ 2\edge[dr] &1\edge[d]
&3\edge[dl] &&&2\edge[ddr] &1\edge[d] &3\edge[ddl] &&&&1 &&&&2\drbl
&1\edge[d] &3\drbl\\
 &1\edge[d] &\dttdar[rrr]\dttdar[uuurrr]\dttdar[dddrrr] &&&&1\edge[d]
&\dttdar[uurrr]\dttdar[ddrrr] &&&& &&&&&1\edge[d]\\
 &1 &&&&&1 &&&&2\drbl &1\edge[d] &3\edge[ddl] &&&&1\\
 &&&&&2\drbl &1\edge[d] &3\edge[dl] &&&&1\edge[d] &\dttdar[uurrr]\\
 &&&&&&1\edge[d] &\dttdar[urrr] &&&&1\\
 &&&&&&1 }$$
\endignore

A justification of this picture is as follows:  One finds
the complement of $\autlap.C$ in the $\grassS$-closure
$\overline{\autlap.C}^\S$  to have the following connected components $U$,
 $V$, $W$:
$$\xalignat2 U &= \{ (k_1,k_2,0,0) \mid k_1, k_2 \in K\}, &V &=  \{ (0,
0, k_3,k_4) \mid k_3, k_4 \in K\}, \\
  W &= \{ (0,k_2,0,k_4) \mid k_2, k_4 \in K\}.
\endxalignat$$  
To further decompose $U$ into
$\autlap$-orbits, we compute the orbit of $(1,0,0,0)$ to be 
$U_1 =
\{k \in U \mid k_1 \in K^*\} \cong K^* \times \AA^1$, leaving the
difference $U \setminus U_1 = (0,k_2,0,0) \mid k_2 \in K\}$.  The latter,
finally, falls into a one-dimensional orbit, $U_2 = \{(0,k_2,0,0) \mid
k_2 \in K^*\}
\cong K^*$, and the point $U_3 = \{(0,0,0,0)\}$, these being the orbits
corresponding to the modules $\bigl( (\la e_1 \oplus \la e_2 ) /
\la(\alpha^2 -
\beta) \bigr) \ \oplus \ S_3$ and $\la e_1 \oplus S_2 \oplus S_3$,
respectively.  The degenerations resulting from this decomposition of
$U$ into orbits are depicted in the uppermost sequence of dotted arrows
in the above diagram.  With
$V$ and $W$ one deals analogously, and on eliminating overlaps, one
obtains the above hierarchy of degenerations with skeleton
$\S = \S_1$.

To obtain the full collection of layer-stable degenerations of $M$
afforded by the closure of $\autlap.C$ in
$\grassSS$, one repeats these computations for the skeleta $\S_2$ and
$\S_3$ of $M$.  The complete hierarchy of proper top-stable degenerations
of $M$ is given in Diagram 5.8 below. The four minimal top-stable
degenerations can be found in the (expanded) left-most column.  The
next column to the right contains the five top-stable degenerations of
height $2$ above $M$, the right-most column the three maximal top-stable
degenerations. \qed
\enddefinition 

\midinsert
\hphantom{}
\medskip
\centerline{\bf The poset of proper top-stable degenerations for Example
5.8}
\ignore
$$\xymatrixrowsep{1.2pc}\xymatrixcolsep{0.1pc}
\xymatrix{
  &&&&&&&& &&&&&&&&&&&&2\edge[d] &&1\drbl &&3\drbl \\
  &&&&&&&& &&&&&&&&&&&&1\edge[d] \\
  &&&&&&&& &&&&&&&&&&&&1 &&&&&&\dttdar[drrrr] \\ 
  &&&&&&&& &&\dttdar[uurrrrrrrr] &&&&&&&&&& &&&&&&&&&& \\
 &&&&2\edge[drr] &&1\edge[d] &&3\drbl  &&&&&&&&&&&&2\edge[ddrr] &&1\edge[d]
&&3\drbl &&&&&&&&2\edge[d] &&1\edge[d] &&3\drbl \\
 &&&& &&1\edge[d]  &&&&\dttdar[rrrrrrrr] &&&&&&&&&&&&1\edge[d] &&&&\dttdar[rrrr]
&&&&&&1 &&1 \\
 &&&& &&1  &&&&&&&&&&&&&&&&1 &&&&\dttdar[drrrr] &&&& \\ 
  &&&&&&\save+<0ex,-3ex> \dttdar[-4,12] \restore && &&&&&&\dttdar[urrrr]
&&&&&&&&&&&&\dttdar[urrrr] &&&& \\ 
2\edge[drr] &&1\drbl &&3\edge[dll]  &&&&2\edge[ddrr] &&1\edge[d] &&3\edge[ddll]
&&&&&&&&2\edge[dr] &&3\edge[dl] &&1\edge[d] &&&&&&&&2\drbl &&1\edge[d] &&3\drbl \\
 &&1\edge[d]  &&&&\save+<0ex,-3ex> \dttdar[6,12] \restore &&&&1\edge[d]
&&&&\dttdar[rrrr] &&&&&&&1 &&&1 &&&&&&&&&&1\edge[d] \\
 &&1  &&&&&&&&1 &&&&\dttdar[drrrr] &&&&&&&& &&&&\dttdar[drrrr]
&&&&&&&&1 \\  
  &&&&&&&& &&&&&& &&&&&&&&&&&&\dttdar[urrrr] &&&& \\
 &&&&2\drbl &&1\edge[d] &&3\edge[dll]  &&&&&&&&&&&&2\drbl &&1\edge[d]
&&3\edge[ddll] &&&&&&&&2\drbl &&1\edge[d] &&3\edge[d] \\
 &&&& &&1\edge[d]  &&&&\dttdar[rrrrrrrr] &&&&&&&&&&&&1\edge[d] &&&&\dttdar[rrrr]
&&&&&&&&1 &&1 \\
 &&&& &&1  &&&&\dttdar[ddrrrrrrrr] &&&&&&&&&&&&1 &&&&&&&& \\ 
  &&&&&&&& &&&&&&&&&&&& &&&&&&\dttdar[urrrr] \\
  &&&&&&&& &&&&&&&&&&&&2\drbl &&1\drbl &&3\edge[d] \\
  &&&&&&&& &&&&&&&&&&&& &&&&1\edge[d] \\
  &&&&&&&& &&&&&&&&&&&& &&&&1 
}$$
\endignore
\medskip
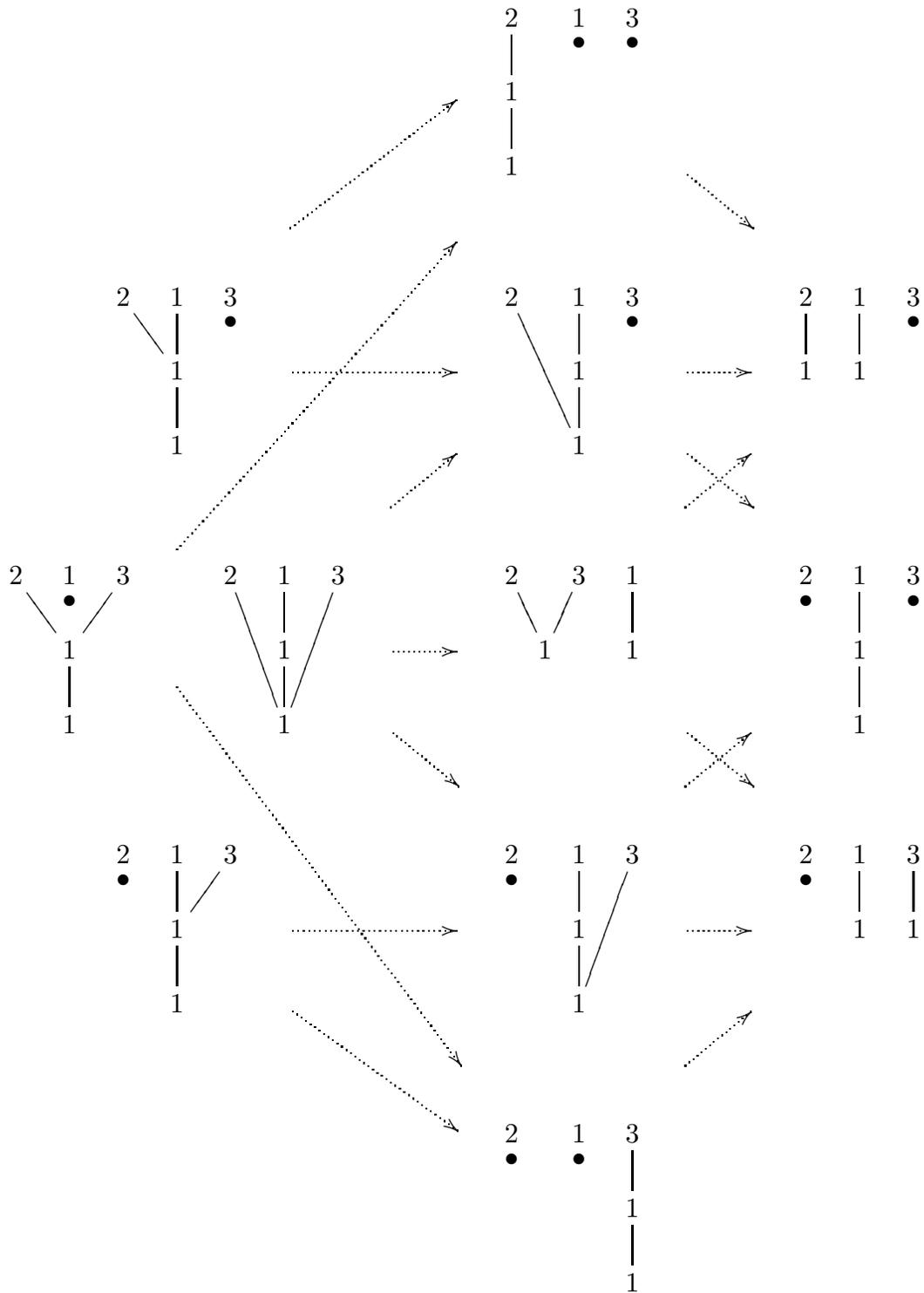
\centerline{Diagram 5.8}
\endinsert

Our final example is a variation on the theme of the previous one, showing
that, for any positive integer $r$, there exists a module $M$ with
squarefree top over a suitable finite dimensional algebra such that
$M$ has a
$\PP^{r-1}$-family of distinct indecomposable layer-stable degenerations.

\definition{Example 5.9} Fix $r\ge 1$, and let $Q$ be the quiver

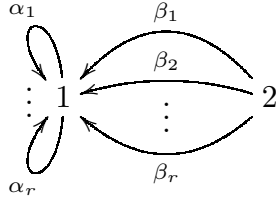
\ignore
$$\xymatrixrowsep{2.0pc}\xymatrixcolsep{2pc}
\xymatrix{
1 \ar@(u,ul)_{\alpha_1} \ar@(d,dl)^{\alpha_r} \save+<-2.5ex,0.5ex>
\drop{\vdots} \restore &\save+<0ex,-0.6ex> \drop{\vdots} \restore
&2 \ar@/_5ex/[ll]_{\beta_1} \ar@/_/[ll]_{\beta_2}
\ar@/^4ex/[ll]^{\beta_r}
}$$
\endignore
  
\noindent Moreover, let $\la =
KQ/I$, where $I$ is the ideal of $KQ$ generated by $\alpha_i^3$ and
$\alpha_i^2\beta_i$ for
$1 \le i \le r$, and $\alpha_i \alpha_j$ and $\alpha_i \beta_j$ for $i \ne
j$.  Moreover, let $P = \la = \la e_1 \oplus \la e_2$ 
and $T$ = $P/JP$.  We focus on the module $M = P/C$ with 
$$C = \sum_{1 \le i
\le r} \la \bigl(\alpha_i - \beta_i \bigr) \in \grasstd,$$
where $d = 2r + 2$.  Here $\autlap.C \cong
\AA^r \times K^*$.

The following is a
$\PP^{r-1}$-family of layer-stable degenerations of $M$:  Namely, for any
point $k \in
\PP^{r-1}$ with homogeneous coordinates $(k_1, \dots, k_r)$, we set $M_k =
P/C_k$, where
$$C_k = \lim_{\tau \rightarrow \infty} C \bigl(e_1 + \tau e_2 + \sum_{1 \le
i \le r} \tau k_i \alpha_i \bigr)= \sum_{1\le i\le r} \la (k_i\alpha_i^2-
\beta_i).$$
Clearly $M_k \cong M_{k'}$ if and
only if $k = k'$ in $\PP^{r-1}$.  For $(k_1, \dots, k_r) \in (K^*)^r$, the
degeneration $M_k$ has the following graph:

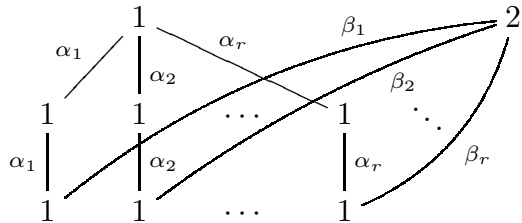
\ignore
$$\xymatrixrowsep{1.5pc}\xymatrixcolsep{1.5pc}
\xymatrix{
 &1 \edge[dl]_{\alpha_1} \edge[d]^(0.6){\alpha_2} \edge[drr]^(0.4){\alpha_r} &&&&2
\edge@/_3ex/[ddlllll]_(0.3){\beta_1} \edge@/_/[ddllll]^(0.3){\beta_2}
\edge@/^2.5ex/[ddll]^{\beta_r}_{\ddots} \\ 
1 \edge[d]_{\alpha_1} &1 \edge[d]^{\alpha_2} &\cdots &1 \edge[d]^{\alpha_r} \\
1 &1 &\cdots &1
}$$
\endignore

\noindent  Beyond the listed ones, $M$ has precisely two further proper
layer-stable degenerations, namely $S_1 \oplus \la e_2$ and $\la e_1 \oplus
S_2$. \qed
\enddefinition

\enddocument